\documentclass{article}
\usepackage[utf8]{inputenc}

\usepackage[T1]{fontenc}
\usepackage{amsmath}
\usepackage{amsfonts}
\usepackage{hyperref}
\hypersetup{
    colorlinks=true,
    linkcolor=cyan,
    filecolor=magenta,      
    urlcolor=cyan,
    pdfpagemode=FullScreen,
    }
\usepackage{amssymb}
\usepackage{mathtools}
\usepackage[version=4]{mhchem}
\usepackage{tabularx}
\usepackage[a4paper, total={6.5in, 9in}]{geometry}
\usepackage[most]{tcolorbox}
\usepackage{xcolor}
\usepackage[numbers]{natbib}
\usepackage[noamsbb]{stmaryrd}
\usepackage{bbold}
\usepackage{multirow} 
\usepackage{enumitem}
\usepackage{tikz-cd}
\usepackage{epstopdf}
\usepackage{bm}
\usepackage{pgfplots}
\usepackage{subcaption}
\usepackage{graphicx}
\usepackage[export]{adjustbox}
\usepackage{amsthm}
\usepackage{chngcntr}
\usepackage{apptools}
\usepackage{esint}
\usepackage{lmodern}
\usepackage{xfrac}
\usepackage{tikz}
\usetikzlibrary{patterns,positioning,cd,calc}
\pgfplotsset{compat=1.18}
\usepackage{authblk}
\usepackage{float}

\newtheorem{theorem}{Theorem}[section]
\newtheorem{lemma}[theorem]{Lemma}

\newtheorem{remark}{Remark}

\newtheorem{assumption}{Assumption}[section]

\tcbset{
    coloredtheoremstyle/.style={
        enhanced,
        breakable,
        sharp corners,
        boxrule=0pt,
        frame hidden,
        colback=black!10!white,
        coltext=black,
        left=6pt,
        right=6pt,
        top=6pt,
        bottom=6pt,
        before skip=10pt,
        after skip=10pt
    }
}
\tcolorboxenvironment{theorem}{coloredtheoremstyle}
\tcolorboxenvironment{proposition}{coloredtheoremstyle}
\tcolorboxenvironment{lemma}{coloredtheoremstyle}
\tcolorboxenvironment{corollary}{coloredtheoremstyle}
\tcolorboxenvironment{claim}{coloredtheoremstyle}
\tcolorboxenvironment{statement}{coloredtheoremstyle}
\tcolorboxenvironment{conjecture}{coloredtheoremstyle}
\tcolorboxenvironment{assumption}{coloredtheoremstyle}
\tcolorboxenvironment{definition}{coloredtheoremstyle}
\newtcolorbox{formulationbox}[1][]{
    enhanced,
    ams equation,
    sharp corners,
    boxrule=0.8pt,
    colback=black!2!white,
    colframe=black!35!white,
    left=6pt,
    right=6pt,
    top=4pt,
    bottom=4pt,
    before skip=10pt,
    after skip=10pt,
    #1
}
\NewDocumentEnvironment{formulation}{O{}}
    {\begin{formulationbox}[#1]}
    {\end{formulationbox}\noindent\ignorespacesafterend}
\numberwithin{equation}{section}

\def\eps{{\varepsilon}}
\DeclareMathOperator*{\esssup}{ess\,sup}

\def\div{\mathrm{div\,}}
\def\divn{\mathrm{div}}
\def\curl{\boldsymbol{\mathrm{curl}}\,}
\def\curln{\boldsymbol{\mathrm{curl}}}

\def\grad{\boldsymbol{\mathrm{grad}}\,}
\def\gradn{\boldsymbol{\mathrm{grad}}}

\def\vecDelta{\boldsymbol{\Delta}}

\def\Hdiv{{\boldsymbol{H}(\mathrm{div})}}
\def\Hzdiv{{\boldsymbol{H}_0(\mathrm{div})}}
\def\Hdivz{{\boldsymbol{H}(\mathrm{div}0)}}
\def\Hzdivz{{\boldsymbol{H}_0(\mathrm{div}0)}}
\def\Hcurl{{\boldsymbol{H}(\boldsymbol{\mathrm{curl}})}}
\def\Hzcurl{\boldsymbol{H}_0(\boldsymbol{\mathrm{curl}})}

\def\ub{{\boldsymbol{u}}}
\def\vb{{\boldsymbol{v}}}

\def\nb{{\boldsymbol{n}}}

\def\R{{\mathbb{R}}}

\def\fb{\boldsymbol{f}}

\def\psib{\boldsymbol{\psi}}
\def\Hb{\boldsymbol{H}}
\def\Lb{\boldsymbol{L}}
\def\asf{\mathsf{a}}

\def\Nedelec{N\'{e}d\'{e}lec }
\def\xb{\boldsymbol{x}}

\def\xib{\boldsymbol{\xi}}
\def\Wb{\boldsymbol{W}}

\def\wb{{\boldsymbol{w}}}

\def\Sigmab{\boldsymbol{\Sigma}}

\def\Vb{\boldsymbol{V}}

\def\Cb{\boldsymbol{C}}

\def\Ecal{\mathcal{E}}

\def\dsf{\mathsf{d}}
\def\csf{\mathsf{c}}
\def\tsf{\mathsf{t}}
\def\Ksf{\mathsf{K}}

\def\ILGav{\mathcal{I}_h}
\def\ILGavz{\mathcal{I}_{h0}}
\def\INDav{\boldsymbol{\mathcal{I}}^{\mathrm{ND}}_h}
\def\INDavz{\boldsymbol{\mathcal{I}}^{\mathrm{ND}}_{h0}}
\def\IRTav{\boldsymbol{\mathcal{I}}^{\mathrm{RT}}_h}
\def\IRTavz{\boldsymbol{\mathcal{I}}^{\mathrm{RT}}_{h0}}
\def\IDGav{\mathcal{I}^{\mathrm{disc}}_h}
\def\IDGavz{\mathcal{I}^{\mathrm{disc}}_{h0}}

\def\Xb{{\boldsymbol{X}}}

\def\Holders{H\"{o}lder's }

\def\zerob{\boldsymbol{0}}
\def\Zb{\boldsymbol{Z}}
\def\Sigmabh{\boldsymbol{\Sigma}_h}

\def\Fcal{\mathcal{F}}

\def\Tcal{\mathcal{T}}

\def\etab{\boldsymbol{\eta}}

\def\Pcal{\mathcal{P}}
\def\Ib{\boldsymbol{I}}

\def\dx{\,\mathrm{d}x}
\def\ds{\,\mathrm{d}s}

\def\Qsf{\mathsf{Q}}

\def\ssf{\mathsf{s}}

\def\Bb{\boldsymbol{B}}
\def\Cb{\boldsymbol{C}}
\def\ab{\boldsymbol{a}}
\def\Rbcal{\boldsymbol{\mathcal{R}}}
\def\eb{\boldsymbol{e}}

\def\Dcal{\mathcal{D}}

\def\element{\mathsf{T}}
\def\face{\mathsf{F}}

\newcommand{\vertiii}[1]{{\left\vert\kern-0.25ex\left\vert\kern-0.25ex\left\vert #1 
    \right\vert\kern-0.25ex\right\vert\kern-0.25ex\right\vert}}
\newcommand{\jump}[1]{\left[\hspace{-1.5pt}\left[#1\right]\hspace{-1.5pt}\right]}
\newcommand{\avg}[1]{\left\{\hspace{-3pt}\left\{#1\right\}\hspace{-3pt}\right\}}

\title{A robust finite element method for incompressible magnetohydrodynamics using Whitney forms}
\author[1]{Tianwei Yu\thanks{\texttt{tianwei.yu@sam.math.ethz.ch}}}
\affil[1]{Seminar for Applied Mathematics, ETH Z\"{u}rich}

\date{}

\begin{document}
\maketitle
\begin{abstract}
A stabilized semi-discrete finite element scheme for the incompressible magnetohydrodynamics system is presented. The scheme uses the Whitney $2$-form and $1$-form for the fluid velocity and magnetic induction, respectively. The key ingredient is a magnetic curl-jump penalty term, which allows us to prove a pre-asymptotic error estimate of order $O(h^{1/2})$ with a constant depending on solution norms but not explicitly on the fluid viscosity and magnetic diffusivity. 
\end{abstract}
\section{Introduction}
Let $\Omega \subset \R^3$ be a bounded connected Lipschitz polyhedral domain with boundary $\Gamma := \partial \Omega$. Let $T > 0$, and set the time interval $I_{\mathrm{T}} := (0,T)$. Set $Q_T := \Omega \times I_{\mathrm{T}}$. Consider the following (scaled) incompressible resistive magnetohydrodynamics (MHD) system~\cite[Eq.~(2.77)]{gerbeau_2006}:
\begin{align}
    \partial_t \ub + (\ub \cdot \gradn) \ub - \nu_\mathrm{v} \vecDelta \ub - (\curl \Bb)\times \Bb  - \grad p&= \fb,  &&\text{in } Q_T, \label{eq:mhd-fluid}\\
    \partial_t \Bb + \nu_\mathrm{m} \curl\curl \Bb - \curl (\ub \times \Bb) &= \zerob, &&\text{in } Q_T, \label{eq:mhd-mag}\\
    \div \ub &= 0, &&\text{in } Q_T, \label{eq:mhd-div-ub}\\
    \div \Bb &= 0, &&\text{in } Q_T, \label{eq:mhd-div-Bb}
\end{align}
subject to the initial and boundary conditions~\cite[Eq.~(2.78)]{gerbeau_2006}:
\begin{align}
    &\ub(\xb,0) = \ub_0(\xb),\qquad\Bb(\xb,0) = \Bb_0(\xb)\; &&\text{on } \Omega\times\{0\}, \label{eq:mhd-ic}\\
    &\ub = \zerob,\qquad \Bb \cdot \nb =0,\qquad \curl\Bb \times \nb = 0, &&\text{on } \Gamma \times I_{\mathrm{T}}.\label{eq:mhd-bc}
\end{align}
Here, $\ub: Q_T \to \R^3$ is the velocity field, $\Bb: Q_T \to \R^3$ is the magnetic induction, $p: Q_T \to \R$ is the pressure, $\fb: Q_T \to \R^3$ is the external force, $\nu_\mathrm{v}, \nu_\mathrm{m} \in (0,1]$ are the dimensionless fluid viscosity and magnetic diffusivity coefficients, respectively. The initial conditions $\ub_0$ and $\Bb_0$ are divergence-free and satisfy the boundary conditions~\eqref{eq:mhd-bc}. The MHD system describes the dynamics of electrically conducting fluids (e.g., plasmas, liquid metals) under the Lorentz force imposed by the electromagnetic field, while the motion of the fluid forms currents that in turn induce a magnetic field; see~\cite{freidberg_2014} for a detailed derivation and physical background. In many practical situations, $\nu_\mathrm{v}$ and $\nu_\mathrm{m}$ (also known as the reciprocals of the fluid and magnetic Reynolds numbers) can be very small~\cite{schaeffer_2017,tzeferacos_2018}. In this regime, convection can dominate diffusion. It is well known that special care is needed for the numerical treatment of such convection-dominated problems (see, e.g.,~\cite{roos_2008}) to alleviate spurious oscillations and improve the accuracy of numerical solutions.

A vast spectrum of numerical methods has been proposed to solve different MHD models (e.g., ideal/resistive, compressible/incompressible, with/without the Hall effect) in pursuit of various features. \emph{Entropy stable} numerical schemes that respect a certain entropy inequality for compressible MHD models have been the subject of sustained research interest; see~\cite{kumar_2012,winters_2016,bohm_2020,wu_2020} and references therein. Recent years have seen a surge of attention in developing \emph{structure-preserving} discretizations which maintain certain physical properties at the discrete level (e.g., vanishing divergence, magnetic/cross helicity conservation). These properties are deemed crucial to get reliable long-time solutions; see~\cite{greif_2010,hu_2017b,hiptmair_2018,hu_2019,gawlik_2022,laakmann_2023,mao_2026} and references therein. Since both the fluid and the electromagnetic field get transported under the fluid flow, discretizations of MHD systems are alternatively approached from the perspective of \emph{advection of differential forms}, extending upwinding techniques from scalar advection problems to differential forms; see, e.g.,~\cite{heumann_2013,pagliantini_2016,wang_2023,luo_2025,li_2026} for advection of differential forms, and~\cite{gerbeau_2000,greif_2010,daveiga_2024,daveiga_2026b} for stabilized schemes for MHD problems. We stress that aforementioned directions are not mutually exclusive; in fact, numerical methods that combine multiple features are often pursued.

In this work, we focus on robustness in convection-dominated regimes and propose a stabilized finite element semi-discretization for the incompressible MHD system~\eqref{eq:mhd-fluid}--\eqref{eq:mhd-bc}. The scheme seeks the discrete velocity and magnetic induction in the Whitney $2$-form (i.e., the lowest-order Raviart--Thomas element) and $1$-form spaces (i.e., the lowest-order \Nedelec element of the first kind), respectively. It is pressure-robust (in the sense of~\cite{john_2017}) and free of a magnetic Lagrange multiplier. In particular, it admits an $O(h^{1/2})$ pre-asymptotic error estimate with $(\nu_\mathrm{v},\nu_\mathrm{m})$-independent constant for exact solutions with relatively low regularity (see Assumption~\ref{ass:reg}). 

It is not new to discretize the velocity and the magnetic induction using $\Hdiv$- and $\Hcurl$-conforming elements, respectively. As the exact divergence-free property of the discrete velocity becomes increasingly appreciated in numerical approximations of Navier--Stokes equations~\cite{john_2018}, inf-sup stable Stokes pairs that respect a discrete Stokes complex are of great interest; see~\cite{neilan_2020} and references therein. Alternatively, one can use the standard $\Hdiv$-conforming Raviart--Thomas or Brezzi--Douglas--Marini element, addressing the $\Hb^1$-continuity by techniques from discontinuous Galerkin (DG) methods~\cite{wang_2007,guzman_2017,barrenechea_2020}. Regarding the approximation of the magnetic induction, the choice of $\Hcurl$-conforming elements is driven by the possible singularity of the electromagnetic field in non-convex domains~\cite{costabel_1990}. With the exact solution failing to be $\Hb^1$-regular, seeking the discrete magnetic induction in $\Hb^1$-conforming element spaces may lead to non-convergence unless treated properly~\cite{costabel_2002,hasler_2004,badia_2012}. An alternative strategy is to approximate the magnetic induction in $\Hcurl$-conforming spaces. Combining this magnetic discretization with stable Stokes pairs, the methods in~\cite{schotzau_2004,prohl_2008} converge to weak solutions in the stationary and transient settings, respectively, without requiring domain convexity. More recent works include~\cite{gao_2019,daveiga_2026b,daveiga_2026}. From a physical viewpoint, the magnetic induction is naturally represented by a $2$-form, which motivates its approximation in an $\Hdiv$-conforming finite element space~\cite{hu_2017b,hu_2019,alonso_2026}. Alternatively, one can discretize, instead of the magnetic induction, a magnetic potential~\cite{hiptmair_2018}.

The highlight of this work is \emph{an $O(h^{1/2})$ pre-asymptotic error estimate with $(\nu_\mathrm{v},\nu_\mathrm{m})$-independent constant for relatively rough solutions of the nonlinear MHD system}. Such a convergence rate is deemed satisfactory as the best general $\Lb^2$-approximation order for the Whitney $1$- and $2$-form spaces used here is $O(h)$ for sufficiently regular fields while half-order loss is commonly encountered in robust estimates for stabilized advection discretizations~\cite[Section~2]{roos_2015}. We are largely motivated by the recent work~\cite{daveiga_2026} whose (third) method achieves a pre-asymptotic convergence rate of $O(h^{k+1/2})$ where $k$ is the complete polynomial degree of the finite element spaces. However, it relies on a conjectured discrete trilinear estimate~\cite[Eq.~(6.3)]{daveiga_2026} (see also~\eqref{eq:daveiga-assumption}), which is not known to hold in general. By a judicious choice of discrete spaces and magnetic stabilization term, we establish a counterpart of the trilinear estimate that suffices for our analysis (see Section~\ref{sec:aux}). We remark that \emph{the success is due to seeking the discrete magnetic induction in the Whitney $1$-form space, while the fluid discretization is rather flexible}. Some properties specific to the low-degree polynomial space are exploited (see Sections~\ref{sec:whitneyprop} and~\ref{sec:aux}). 

Lastly, we comment on different settings for MHD error analysis. This work analyzes a fully nonlinear semi-discretization; see also~\cite{daveiga_2025,daveiga_2026}. Conceptually speaking, the presence of the time-derivative terms provides sufficient $\Lb^2$-control thanks to Gronwall's inequality. Without the time-derivatives, error estimates of stationary problems are derived either under small data conditions (see, e.g.,~\cite{schotzau_2004,greif_2010,hu_2020}) or for linearized systems arising from Picard iterations (see, e.g.,~\cite{gerbeau_2000,codina_2006,daveiga_2024}). 

The remainder of this work is organized as follows. We conclude this section by introducing the notation used throughout the paper. In Section~\ref{sec:continuous}, the weak formulation of the MHD system~\eqref{eq:mhd-fluid}--\eqref{eq:mhd-bc} is briefly presented and the regularity assumptions necessary for our analysis are stated. 
Section~\ref{sec:prelim} is devoted to a recap of some machinery of \emph{finite element exterior calculus} (FEEC) and the properties of Whitney forms. In Section~\ref{sec:femdisc}, we present the proposed semi-discrete scheme, and derive the main a priori error estimate in Section~\ref{sec:err-est}. Numerical experiments are presented in Section~\ref{sec:numexp} with concluding remarks in Section~\ref{sec:conclusion}. 

\subsection{Notation}
Given a generic open bounded Lipschitz domain $D \subset \R^3$, we use the standard notation for Sobolev spaces $W^{s,p}(D), s\in \R, p \in [1,\infty]$ and write $H^s(D) := W^{s,2}(D)$ and $L^p(D) := W^{0,p}(D)$. These spaces are equipped with the standard norms $\|\cdot\|_{W^{s,p}(D)}$ and semi-norms $|\cdot|_{W^{s,p}(D)}$. We denote $\Hb(\curln, D)$ and $\Hb(\divn, D)$ as the spaces of square-integrable vector fields with square-integrable curl and divergence, respectively. We also write $H^1_0(D), \Hb_0(\curln, D)$ and $\Hb_0(\divn, D)$ as the subspaces with vanishing value, tangential component, and normal component on $\partial D$, respectively. These conditions are understood in the sense of traces; see, e.g.~\cite[Section~3.4]{arnold_2018}. Additionally, we denote $L^2_0(D)$ as the subspace of $L^2(D)$ with vanishing mean value. The vector versions of the above spaces are denoted by boldface letters, e.g., $\Lb^2(D) := [L^2(D)]^3$. We write $(\cdot,\cdot)_D$ for the $L^2$-inner product on $D$ and $\langle\cdot,\cdot\rangle_{\face}$ for the $L^2$-inner product on a two-dimensional face $\face$. Let $\Tcal = \{\element\}$, and $\Fcal = \{\face\}$ be sets of subdomains $\element$ and faces $\face$. We denote $(u,v)_{\Tcal} := \sum_{\element \in \Tcal} (u,v)_\element$ and $\langle u,v\rangle_{\Fcal} := \sum_{\face \in \Fcal} \langle u,v\rangle_\face$. The broken $L^2$-norms are defined by $\|u\|^2_{L^2(\Tcal)} := \sum_{\element \in \Tcal} \|u\|^2_\element$ and $\|u\|^2_{L^2(\Fcal)} := \sum_{\face \in \Fcal} \|u\|^2_\face$. In all the above notation, we omit the domain $D$ when $D = \Omega$. 

The symbol $C$ stands for a generic positive constant that may change in each occurrence. It may depend on the spatial domain $\Omega$, the mesh shape-regularity, and so on, but always independent of the mesh size $h$, the fluid viscosity $\nu_\mathrm{v}$, and the magnetic diffusivity $\nu_\mathrm{m}$. 

For a Banach space $X$, the Bochner space $L^p(I_{\mathrm{T}};X)$ contains strongly measurable functions with finite $\int_0^T \|u(t)\|_X^p\ \dsf t$ for $1 \leq p < \infty$, and $\esssup_{t\in I_{\mathrm{T}}} \|u(t)\|_X$ for $p = \infty$. Let $H^1(I_{\mathrm{T}};X) := \{u \in L^2(I_{\mathrm{T}};X): \partial_t u \in L^2(I_{\mathrm{T}};X)\}$ where $\partial_t u$ is the weak time derivative. See, e.g.,~\cite[Chapter~64]{ern_2021c} for more details.

For the reader's convenience, some relevant vector calculus identities are listed below:
\begin{align}
    (\xib \times \etab) \cdot \psib &= - (\xib \times \psib) \cdot \etab, \label{eq:triple}\\
    \curl (a \xib) &= \grad a \times \xib + a \curl \xib, \label{eq:curlau}\\
    \curl (\xib \times \etab) &= (\div \etab)\xib - (\div \xib)\etab + (\etab\cdot\gradn)\xib - (\xib\cdot\gradn)\etab, \label{eq:curlutv}\\
    \div(\xib\otimes\etab) &= (\etab \cdot \gradn) \xib + \xib \div \etab, \label{eq:divutv} \\
    \div(\xib\times \etab) &= \etab \cdot \curl \xib - \xib \cdot \curl \etab. \label{eq:divutimesv}
\end{align}
Here, $[\xib \otimes \etab]_{i,j} := \xib_i \etab_j$ denotes the tensor product, and $[(\etab \cdot \gradn)\xib]_i := \sum_{j=1}^3\etab_j \partial_j \xib_i$ is the directional derivative of $\xib$ along $\etab$ componentwise. 

Because we consider only spatial semi-discretization, we suppress the time dependence of most quantities. For instance, $\xib \in \Xb$ is understood as $\xib(\cdot,t) \in \Xb$ for almost every $t \in I_{\mathrm{T}}$.

\section{Continuous problem}\label{sec:continuous}
Set the function spaces:
\begin{equation}
    \Zb^\mathrm{v} := \Hb^1_0 \cap \Hzdivz,\qquad \Zb^\mathrm{m} := \Hcurl \cap \Hzdivz,
\end{equation}
where $\Hzdivz := \{\vb \in \Hzdiv : \div \vb = 0\}$.
Define the trilinear forms:
\begin{align}
    \tsf(\psib; \xib, \etab) &:= ((\psib \cdot \gradn) \xib, \etab), \label{eq:tsf}\\
    \csf(\psib; \xib, \etab) &:= ((\curl \psib) \times \xib, \etab). \label{eq:csf}
\end{align}
The former is related to the fluid convection term $(\ub \cdot \gradn) \ub$ while the latter is related to the Lorentz force $\curl \Bb \times \Bb$ and the magnetic convection term $\curl (\ub \times \Bb)$. 
The following lemma recalls the boundedness and skew-symmetry of $\tsf$ and $\csf$.
\begin{lemma}[boundedness and skew-symmetry]\label{lem:tsf-csf}
    \phantom{aa}
    \begin{enumerate}[label=$(\roman*)$]
        \item Let $\psib, \xib, \etab \in \Hb^1$. There exists a positive constant $C$ such that
        \begin{equation}\label{eq:tsf-bdd}
            |\tsf(\psib; \xib, \etab)| \leq C \|\psib\|_{\Hb^1} \|\xib\|_{\Hb^1} \|\etab\|_{\Hb^1}.
        \end{equation}
        If $\psib \in \Hzdivz$ in addition, then
        \begin{equation}\label{eq:tsf-skew}
            \tsf(\psib; \xib, \etab) = - \tsf(\psib; \etab, \xib).
        \end{equation}

        \item Let $\psib \in \Hcurl$, $\xib \in \Zb^\mathrm{m}$ and $\etab \in \Hb^1$. There exists a positive constant $C$ such that
        \begin{equation}\label{eq:csf-bdd}
            |\csf(\psib; \xib, \etab)| \leq C \|\psib\|_{\Hcurl} \|\xib\|_{\Hcurl} \|\etab\|_{\Hb^1}.
        \end{equation}
        In addition, we have
        \begin{equation}\label{eq:csf-skew}
            \csf(\psib; \xib, \etab) = -\csf(\psib; \etab, \xib).
        \end{equation}
    \end{enumerate}
\end{lemma}
\begin{proof}
    Both boundedness estimates~\eqref{eq:tsf-bdd} and~\eqref{eq:csf-bdd} follow from \Holders inequality, Sobolev embedding $H^1 \hookrightarrow L^6$ (see, e.g.,~\cite[Section~9.3]{brezis_2011}), and the regularity pickup $\Hcurl\cap\Hzdivz \hookrightarrow \Hb^{1/2+\delta}, \delta > 0$~\cite{amrouche_1998}; see details in~\cite[Lemmas~2.5 and 2.6]{schotzau_2004}. The skew-symmetry~\eqref{eq:csf-skew} is a direct consequence of the vector calculus identity~\eqref{eq:triple}; the skew-symmetry~\eqref{eq:tsf-skew} follows from the identity~\eqref{eq:divutv} and integration by parts, along with the fact that $\psib\in\Hzdivz$.
\end{proof}

Let the external force $\fb \in L^\infty(I_{\mathrm{T}}; \Lb^2)$. Following~\cite[Eqs.~(2.88)--(2.90)]{gerbeau_2006}, the weak formulation of the MHD system~\eqref{eq:mhd-fluid}--\eqref{eq:mhd-bc} reads: 
\begin{formulation}[top=-0.5em]\label{eq:MHD-weak}
    \begin{aligned}
        \hspace{-1.5em}\text{Seek }\ub \in L^\infty(I_{\mathrm{T}}; \Hzdivz) \cap L^2(I_{\mathrm{T}}; \Zb^\mathrm{v}),\;\Bb \in L^\infty(I_{\mathrm{T}}; \Hzdivz) &\cap L^2(I_{\mathrm{T}}; \Zb^\mathrm{m}) \text{ such that} \\
        \langle\partial_t \ub, \vb\rangle_{(\Zb^\mathrm{v})^\prime,\Zb^\mathrm{v}} + \nu_\mathrm{v} (\grad \ub, \grad \vb) + \tsf(\ub; \ub, \vb) - \csf(\Bb; \Bb, \vb) &= (\fb, \vb), &&\hspace{-4em}\forall\,\vb \in \Zb^\mathrm{v},\\
        \langle\partial_t \Bb, \Cb\rangle_{(\Zb^\mathrm{m})^\prime,\Zb^\mathrm{m}} + \nu_\mathrm{m} (\curl\Bb, \curl\Cb) + \csf(\Cb; \Bb, \ub)  &= 0, &&\hspace{-4em}\forall\,\Cb \in \Zb^\mathrm{m},\\
        \text{in the sense of distributions in time, with initial conditions } \ub(\cdot,0) &= \ub_0,\; \Bb(\cdot,0) = \Bb_0.
    \end{aligned}
\end{formulation}
Note that the pressure $p$ is eliminated by restricting the velocity space to the divergence-free subspace. To make sense of the duality pairings and the initial conditions, we refer to~\cite[Remark~2.1.5]{gerbeau_2006}. See~\cite[Section~2.2]{gerbeau_2006} for discussions on existence and uniqueness of the weak solution to the MHD system~\eqref{eq:mhd-fluid}--\eqref{eq:mhd-bc}, and also~\cite{schotzau_2004} for a related result in the stationary case.

 We finalize this section by stating the regularity assumptions on the exact solutions used in the error analysis in Section~\ref{sec:err-est}.
\begin{assumption}[regularity of the exact solution]\label{ass:reg}
    For some constant $\delta > 0$, the exact solution $(\ub, \Bb)$ to the MHD system~\eqref{eq:mhd-fluid}--\eqref{eq:mhd-bc} satisfies
    \begin{enumerate}[label=$(\roman*)$]
    \item $\ub \in H^1(I_{\mathrm{T}}; \Hb^{1/2+\delta})\cap L^2(I_{\mathrm{T}}; \Hb^{3/2+\delta}),\; \Bb \in H^1(I_{\mathrm{T}}; \Hb^{1/2+\delta}),\;\curl\Bb \in L^2(I_{\mathrm{T}}; \Hb^{1/2+\delta})$;\label{ass:reg-1}
    \item $\grad\ub$, $\Bb$, $\curl\Bb \in L^\infty(I_{\mathrm{T}}; \Lb^{\infty})$.\label{ass:reg-2}
    \end{enumerate}
\end{assumption}
\begin{remark}[on the regularity assumption]\label{rmk:reg}\phantom{aa}
    The spatial regularities in Assumption~\ref{ass:reg}\,\ref{ass:reg-1} are, despite the lack of a full proof, deemed minimal provided the external force $\fb \in L^\infty(I_{\mathrm{T}}; \Lb^2)$ on polyhedral domains; see discussion in~\cite[Remark~4.1]{greif_2010}. We impose additionally the $H^1$ time regularity in~\ref{ass:reg}\,\ref{ass:reg-1} and $\Lb^\infty$ space-time regularity in Assumption~\ref{ass:reg}\,\ref{ass:reg-2} for our convenience in the error analysis~\eqref{eq:err-est}. We remark that it is possible to relax them. 
\end{remark}

\section{Preliminaries}\label{sec:prelim}
In this section, we first introduce the finite element spaces, followed by some well-known interpolation estimates. We shall also present some auxiliary results that are crucial for the error analysis in Section~\ref{sec:err-est}.

\subsection{Finite element differential forms}
Let $\{\Tcal_h\}_{0<h\leq 1}$ be a family of shape-regular and quasi-uniform simplicial triangulations of $\Omega$ with mesh size $h$. We denote by $\Fcal_h$ the set of all faces in $\Tcal_h$ and by $\Fcal^\partial_h$ and $\Fcal^\circ_h$ the sets of boundary and interior faces, respectively. Each face $\face \in \Fcal_h$ is endowed with a unit normal vector $\nb_\face$ and, in particular, $\nb_\face = \nb$ for $\face \in \Fcal^\partial_h$. Given an interior face $\face \in \Fcal^\circ_h$, we denote by $\element^+$ and $\element^-$ the two elements sharing $\face$ and such that  $\nb_\face$ is pointing from $\element^+$ to $\element^-$. 
Denote $v^+ := v|_{\element^+}$ and $v^- := v|_{\element^-}$. We introduce the notion of jump and average operators for a piecewise smooth function $v$ as:
\begin{equation}
    \jump{v}_\face := v^+|_\face - v^-|_\face,\qquad \avg{v}_\face := \frac{1}{2}\left(v^+|_\face + v^-|_\face\right).
\end{equation}
For a boundary face $\face \in \Fcal^\partial_h$, we adopt the convention that $\jump{v}_\face = \avg{v}_\face = v|_\face$.

Throughout the manuscript, we adopt the following notation for the Whitney forms~\cite{arnold_2018}: 
\begin{equation}\label{def:whitney}
    \begin{aligned}
        P_h &:= \Pcal_1^-\Lambda^0 :=  \{v \in H^1:v|_\element \in \mathcal{P}_1\;\;\forall\,\element\in\mathcal{T}_h\}, & P_{h,0} &:= P_h\cap H^1_0, \\
        \Sigmabh &:= \Pcal_1^-\Lambda^1 := \{\vb \in \Hcurl: \vb|_\element \in \boldsymbol{\mathcal{P}}_0 + \xb\times\boldsymbol{\mathcal{P}}_0\;\;\forall\,\element\in\mathcal{T}_h\}, & \Sigmab_{h,0} &:= \Sigmab_{h}\cap\Hzcurl,\\
        \Vb_{h} &:= \Pcal_1^-\Lambda^2 := \{\vb \in \Hdiv: \vb|_\element \in \boldsymbol{\mathcal{P}}_0 + \xb\mathcal{P}_0\;\;\forall\,\element\in\mathcal{T}_h\},  & \Vb_{h,0} &:= \Vb_h\cap\Hzdiv, \\
        S_h &:= \Pcal_1^-\Lambda^3 := \{v \in L^2: v|_\element \in \mathcal{P}_0\;\;\forall\,\element\in\mathcal{T}_h\}, & S_{h,0} &:= S_h\cap L^2_0,
    \end{aligned}
\end{equation}
where $\mathcal{P}_k$ denotes the space of complete 3-variate polynomials with total degree at most $k$.

In the proposed scheme~\eqref{eq:semidiscrete}, the discrete velocity $\ub_h$ and magnetic induction $\Bb_h$ are sought in the spaces $\Vb_{h,0}$ and $\Sigmab_{h}$, respectively. To accommodate the divergence-free constraints~\eqref{eq:mhd-div-ub} and~\eqref{eq:mhd-div-Bb} at the discrete level, we define the following subspaces:
\begin{equation}
    \Zb^{\mathrm{v}}_h := \Vb_{h,0} \cap \Hdivz,\qquad \Zb^{\mathrm{m}}_h := \Sigmab_{h} \cap (\grad P_h)^\perp
\end{equation}
where $\perp$ stands for the $\Lb^2$-orthogonal complement. Note the non-conformity  $\Zb^{\mathrm{v}}_h \not\subset \Zb^{\mathrm{v}}$ and $\Zb^{\mathrm{m}}_h \not\subset \Zb^{\mathrm{m}}$.

\subsection{Bounded cochain projections}\label{sec:cochainproj}
We recall the de Rham complex and the existence of well-established $\Lb^2$-bounded \emph{cochain projections} (and their variants with zero trace) in the following commuting diagrams~\cite[Theorem~23.12]{ern_2021}:
\begin{equation*}
    \hspace{-3em}
    \begin{tikzcd}[row sep=1.5em, column sep=3em]
        & H^1 \arrow[d, "\ILGav"] \arrow[r, "\grad"]
        & \Hcurl \arrow[d, "\INDav"] \arrow[r, "\curl"]
        & \Hdiv \arrow[d, "\IRTav"] \arrow[r, "\div"]
        & L^2 \arrow[d, "\IDGav"]\\
        & P_{h} \arrow[r, "\grad"]
        & \Sigmab_{h} \arrow[r, "\curl"]
        & \Vb_{h} \arrow[r, "\div"]
        & S_{h} 
    \end{tikzcd}
    \hspace{-3em}
    \begin{tikzcd}[row sep=1.5em, column sep=3em]
        & H^1_0 \arrow[d, "\ILGavz"] \arrow[r, "\grad"]
        & \Hzcurl \arrow[d, "\INDavz"] \arrow[r, "\curl"]
        & \Hzdiv \arrow[d, "\IRTavz"] \arrow[r, "\div"]
        & L^2_0 \arrow[d, "\IDGavz"]\\
        & P_{h,0} \arrow[r, "\grad"]
        & \Sigmab_{h,0} \arrow[r, "\curl"]
        & \Vb_{h,0} \arrow[r, "\div"]
        & S_{h,0} 
    \end{tikzcd}
\end{equation*}
We recall some relevant approximation properties (due to~\cite[Theorem~23.12, Corollary~22.9 and Corollary~22.16]{ern_2021}). The following estimates hold whenever the right-hand side is finite: 
\begin{align}
    \qquad\qquad\|\vb - \INDav \vb\|_{\Lb^2} &\leq Ch^{s}|\vb|_{\Hb^{s}} && s\in [0,1],\qquad\qquad \label{eq:approxl2ND}\\
    \qquad\qquad\|\vb - \IRTav \vb\|_{\Lb^2} &\leq Ch^{s}|\vb|_{\Hb^{s}} && s\in [0,1],\qquad\qquad \label{eq:approxl2RT}\\
    \qquad\qquad\|\vb - \IRTav \vb\|_{\Lb^\infty} &\leq Ch^s |\vb|_{\Wb^{s,\infty}} && s \in \{0,1\},\qquad\qquad \label{eq:approxlinfRT}
\end{align}
For the zero-trace cochain projection $\IRTavz$, it holds that, when $\vb \in \Hzdiv$, 
\begin{align}
    \qquad\qquad\|\vb - \IRTavz \vb\|_{\Lb^2} &\leq Ch^s |\vb|_{\Hb^{s}} && s \in [0, 1],\qquad\qquad \label{eq:approxl2RT0}\\
    \qquad\qquad\|\vb - \IRTavz \vb\|_{\Lb^\infty} &\leq Ch^s |\vb|_{\Wb^{s,\infty}} && s \in \{0,1\}.\qquad\qquad \label{eq:approxlinfRT0}
\end{align}
Combining the above estimates with the commuting properties of the above diagrams yields
\begin{align}
    \qquad\qquad\|\curl (\vb - \INDav \vb)\|_{\Lb^2} &\leq Ch^s|\curl \vb|_{\Hb^{s}}&& s\in [0,1],\qquad\qquad \label{eq:approxcurlND}\\
    \qquad\qquad\|\curl (\vb - \INDav \vb)\|_{\Lb^\infty} &\leq Ch^s|\curl \vb|_{\Wb^{s,\infty}}&& s\in \{0,1\}.\qquad\qquad \label{eq:approxcurlinfND}
\end{align}

\subsection{An orthogonal interpolation}\label{sec:rh}
The magnetic induction $\Bb$ lives in $\Zb^\mathrm{m}$ while, as can be seen later in~\eqref{eq:semidiscrete}, the discrete magnetic induction $\Bb_h$ stays in $\Zb_h^{\mathrm{m}}$. Therefore, we shall need an interpolation operator that maps $\Zb^\mathrm{m}$ to $\Zb_h^{\mathrm{m}}$. The cochain projection $\INDav$ in the previous section does not in general serve this purpose.

Let $\Qsf_h: \Lb^2 \rightarrow \grad P_h$ be the $\Lb^2$-orthogonal projection. We define the interpolation operator $\Rbcal^\mathrm{ND}_h:\Lb^2 \rightarrow \Sigmab_h$ as
\begin{equation}
    \Rbcal^\mathrm{ND}_h \vb := \INDav \vb + \Qsf_h (\vb - \INDav \vb).
\end{equation}
It is easy to see the orthogonality property (which implies $\Rbcal^\mathrm{ND}_h (\Zb^\mathrm{m}) \subset \Zb_h^{\mathrm{m}}$):
\begin{equation}\label{eq:orthgradR}
    (\vb - \Rbcal^\mathrm{ND}_h \vb, \grad w_h) = 0 \qquad \forall\,w_h \in P_h, 
\end{equation}
the commuting property:
\begin{equation}\label{eq:commrbcal}
    \curl \Rbcal^\mathrm{ND}_h \vb = \IRTav \curl \vb \qquad \forall\,\vb \in \Hcurl.
\end{equation}
Consequently, $\Rbcal^\mathrm{ND}_h$ inherits the approximation properties of $\INDav$ from Section~\ref{sec:cochainproj}, satisfying the following estimates whenever the right-hand side is finite:
\begin{align}
    \qquad\qquad\|\vb - \Rbcal^\mathrm{ND}_h \vb\|_{\Lb^2} &\leq Ch^{s}|\vb|_{\Hb^{s}} && s\in [0,1],\qquad\qquad \label{eq:approxl2R}\\
    \qquad\qquad\|\curl (\vb - \Rbcal^\mathrm{ND}_h \vb)\|_{\Lb^2} &\leq Ch^s|\curl \vb|_{\Hb^{s}}&& s\in [0,1],\qquad\qquad \label{eq:approxcurlR}\\
    \qquad\qquad\|\curl (\vb - \Rbcal^\mathrm{ND}_h \vb)\|_{\Lb^\infty} &\leq Ch^s|\curl \vb|_{\Wb^{s,\infty}}&& s\in \{0,1\}.\qquad\qquad \label{eq:approxcurlinfR}
\end{align}

\subsection{Some properties of the Whitney forms}\label{sec:whitneyprop}
We collect properties of low-order differential-form spaces that are useful in the analysis below.

The first one is an elementary fact about the Whitney $2$-form $\Pcal_1^-\Lambda^2$.
\begin{lemma}[vanishing elementwise gradient]
    \begin{equation}\label{eq:RTzerograd}
        [\;\vb_h \in \Zb^{\mathrm{v}}_h\;] \implies [\;\grad \vb_h|_\element = 0\quad \forall\,\element \in \Tcal_h\;].
    \end{equation}
\end{lemma}
\begin{proof}
    By definition of $\Vb_h$, we can write $\vb_h|_\element = \ab_\element + b_\element\xb$ for some $\ab_\element \in \R^3$ and $b_\element \in \R$. By the identity~\eqref{eq:curlau}, we have $\curl\vb_h|_\element = 0$. Since $\div \vb_h|_\element = 3b_\element$ and $\grad \vb_h|_\element = b_\element \Ib$ (where $\Ib \in \R^{3 \times 3}$ is the identity matrix), $\vb_h \in \Zb^{\mathrm{v}}_h$ implies $b_\element = 0$ and thus $\grad \vb_h|_\element = 0$. 
\end{proof}

The next lemma bounds the curl seminorm using the $\Lb^2$-norm and tangential jumps of the curl. This estimate is central to the magnetic stabilization and error analysis below.
\begin{lemma}[discrete curl inverse inequality]\label{lem:curlinv}
    There exists a constant $C > 0$ independent of $h$ such that
    \begin{equation}\label{eq:curlinv}
        \|\curl \vb_h\|^2_{\Lb^2} \leq Ch^{-1/2}\|\vb_h\|_{\Lb^2}\left\|\jump{\curl\vb_h}\times\nb\right\|_{\Lb^2(\Fcal_h)} \qquad \forall\,\vb_h \in \Sigmab_h.
    \end{equation}
\end{lemma}
\begin{proof}
    The estimate is a consequence of integration by parts and the lowest polynomial degree. Indeed, for each element $\element \in \Tcal_h$, we have
    \begin{equation}
        \begin{aligned}
            \|\curl \vb_h\|^2_{\Lb^2(\element)} = \int_\element (\curl \vb_h)\cdot(\curl \vb_h)\dx = \int_\element (\curl \curl \vb_h) \cdot \vb_h \dx - \int_{\partial \element}  (\vb_h\times \nb_\element) \cdot \curl \vb_h  \ds
        \end{aligned}
    \end{equation}
    where $\nb_\element$ is the unit outward normal vector on $\partial \element$. 
    The first term vanishes. Summing over all elements $\element \in \Tcal_h$ yields
    \begin{equation}
        \|\curl\vb_h\|_{\Lb^2}^2 = - \sum_{\face \in \Fcal^\circ_h} \int_\face \jump{\vb_h}\times\nb_\face \cdot \avg{\curl\vb_h}_\face + \avg{\vb_h\times\nb}_\face \cdot \jump{\curl\vb_h}_\face  \ds - \int_{\Gamma} (\vb_h\times \nb) \cdot \curl \vb_h \ds.
    \end{equation}
    The interior face jump $\jump{\vb_h}\times\nb_\face$ for $\face \in \Fcal^\circ_h$ vanishes since $\vb_h \in \Hcurl$. The desired result follows from the trace inverse inequality.
\end{proof}
\noindent In fact, we remark that the above result holds for $\Pcal_1 \Lambda^1$ as well.

\section{Finite element discretization}\label{sec:femdisc}
This section presents the proposed finite element semi-discretization of the MHD system~\eqref{eq:mhd-fluid}--\eqref{eq:mhd-bc} and conducts an a priori error analysis. The following paragraphs motivate the discretization of each term using standard techniques. The complete scheme is given in Section~\ref{sec:scheme}.

\subsection{Discretization of individual terms}\label{sec:indivterms}
\paragraph{Fluid diffusion.}
Due to the non-conformity $\Vb_{h,0} \not\subset \Hb^1$, the fluid diffusion term $\nu_\mathrm{v} \vecDelta \ub$ (combined with the pressure term $\grad p$) is treated in a DG fashion. Although this approach is well-established (see, e.g.,~\cite{toselli_2002,wang_2007,greif_2010}), a brief derivation is provided for completeness.
By the regularity stated in Assumption~\ref{ass:reg}, it holds that $\nu_\mathrm{v} \grad \ub + p\Ib \in \Hdiv$ with $\Ib\in\R^{3\times3}$ being the identity matrix. The regularity theory of the Stokes problem~\cite[Theorem~9.20(a)]{dauge_1989} further implies that $\ub \in \Hb^{3/2+\eps}$ and $p \in H^{1/2+\eps}$. Thus, we can perform integration by parts elementwise for $\vb_h \in \Vb_{h,0}$ and write
\begin{equation}\label{eq:intmdeltauv}
    (- \div(\nu_\mathrm{v} \grad \ub + p\Ib), \vb_h)_{\Tcal_h} = (\nu_\mathrm{v} \grad \ub + p\Ib, \grad \vb_h)_{\Tcal_h} - \sum_{\element \in \Tcal_h} \langle (\nu_\mathrm{v} \grad \ub + p\Ib)\nb_\element, \vb_h\rangle_{\partial \element}
\end{equation}
where $\nb_\element$ is the unit outward normal vector on $\partial \element$. Using the jump and average operators, the second term on the right-hand side of~\eqref{eq:intmdeltauv} evaluates to
\begin{equation}\label{eq:stokes-jumpterm}
    \begin{aligned}
        -\langle \jump{\nu_\mathrm{v}\grad \ub + p\Ib}\nb, \avg{\vb_h}\rangle_{\Fcal^\circ_h} - \langle \avg{\nu_\mathrm{v}\grad \ub + p\Ib}\nb, \jump{\vb_h}\rangle_{\Fcal^\circ_h} - \langle (\nu_\mathrm{v}\grad \ub + p\Ib)\nb, \vb_h\rangle_{\Gamma}.
    \end{aligned}
\end{equation}
The first term in~\eqref{eq:stokes-jumpterm} vanishes. Also, $\langle \avg{p\Ib}_\face\nb_\face, \jump{\vb_h}_\face\rangle_\face = 0$ for $\face \in \Fcal_h$ since $\vb_h \in \Vb_{h,0} \subset \Hzdiv$. Consequently,~\eqref{eq:intmdeltauv} reduces to
\begin{equation}\label{eq:mdeltauv}
        (- \div(\nu_\mathrm{v} \grad \ub + p\Ib), \vb_h)_{\Tcal_h} = \nu_\mathrm{v}(\grad \ub, \grad \vb_h)_{\Tcal_h} + (p, \div \vb_h) - \nu_\mathrm{v}\langle \avg{\grad \ub}\nb, \jump{\vb_h}\rangle_{\Fcal_h}.
\end{equation}
We adopt the \emph{symmetric interior penalty} (SIP) (see, e.g.,~\cite[Section~4.2]{pietro_2012}) in the standard DG formulation. That is, a symmetric term and a penalty term are added to~\eqref{eq:mdeltauv}. Define the bilinear form:
\begin{equation}\label{eq:asfh}
    \asf_h(\xib, \etab) := (\grad \xib, \grad \etab)_{\Tcal_h} - \langle \avg{\grad \xib}\nb, \jump{\etab}\rangle_{\Fcal_h} - \langle \avg{\grad \etab}\nb, \jump{\xib}\rangle_{\Fcal_h} + \gamma_\asf h^{-1}\langle \jump{\xib}, \jump{\etab}\rangle_{\Fcal_h}
\end{equation}
where $\gamma_\asf > 0$ is a penalty parameter. Combining the above arguments, for the exact $\ub$ and $p$, we obtain the consistency:
\begin{equation}\label{eq:vecdiffusion-consistency}
    (- \div(\nu_\mathrm{v} \grad \ub + p\Ib), \vb_h)_{\Tcal_h} = \nu_\mathrm{v}\asf_h(\ub, \vb_h) + (p, \div \vb_h) \qquad \forall\,\vb_h \in \Vb_{h,0}.
\end{equation}
It is standard to check the coercivity and boundedness of $\asf_h$. 
\begin{lemma}[coercivity and boundedness of $\asf_h$]\label{lem:asf}
    There exist constants $\gamma^\ast_\asf,\beta_\asf, C_\asf > 0$ independent of $h$ such that if $\gamma_\asf > \gamma^\ast_\asf$ in~\eqref{eq:asfh}, it holds that
    \begin{align}
        \qquad\qquad\asf_h(\vb_h, \vb_h) &\geq \beta_\asf \|\vb_h\|^2_{\asf,h} &&\forall\,\vb_h \in \Vb_{h,0},\qquad\qquad \\
        \qquad\qquad|\asf_h(\xib, \etab)| &\leq C_\asf \vertiii{\xib}_{\asf,h} \vertiii{\etab}_{\asf,h} && \forall\,\xib, \etab \in \Hb^{3/2+\delta}(\Tcal_h)+\Vb_{h,0},\qquad\qquad
    \end{align}
    where the norms are defined as
    \begin{equation}\label{eq:asf-norm}
        \|\xib\|^2_{\asf,h} := \|\grad \xib\|^2_{\Tcal_h} + h^{-1}\|\jump{\xib}\|^2_{\Fcal_h},\qquad \vertiii{\xib}^2_{\asf,h} := \|\xib\|^2_{\asf,h} + h\|\avg{\grad \xib}\|^2_{\Fcal_h}.
    \end{equation}
\end{lemma}
\begin{proof}
    The proof is elementary; see, e.g.,~\cite[Lemmas~4.1 and~4.2]{wang_2007}.
\end{proof}
\noindent Note that a sufficiently large $\gamma_\asf > \gamma^\ast_\asf$ as in Lemma~\ref{lem:asf} is taken for granted in the subsequent analysis.

\paragraph{Fluid convection.}
Due to the non-conformity $\Vb_{h,0} \not\subset \Hb^1$, the trilinear form $\tsf$ in~\eqref{eq:tsf} is not well-defined on $\Vb_{h,0}\times\Vb_{h,0}\times\Vb_{h,0}$. 
Define the following discrete trilinear form:
\begin{equation}\label{eq:tsfh}
    \tsf_h(\psib; \xib, \etab) := ((\psib \cdot \gradn) \xib, \etab)_{\Tcal_h} - \langle (\psib \cdot \nb) \jump{\xib}, \avg{\etab}\rangle_{\Fcal^\circ_h}.
\end{equation}
The interior face integral term is motivated by the fact that, when $\psib \in \Hzdivz$, the trilinear form $\tsf_h$ is skew-symmetric in the last two arguments, mimicking the continuous case~\eqref{eq:tsf-skew} (see, e.g.,~\cite[Lemma~2.1]{guzman_2017} in the context of incompressible Euler's equations). This is stated in the following lemma.
\begin{lemma}[skew-symmetry of $\tsf_h$]
    Let $\psib \in \Hb^1(\Tcal_h)\cap\Hzdivz$. It holds that
    \begin{equation}\label{eq:tsfh-skew}
        \tsf_h(\psib; \xib, \etab) = - \tsf_h(\psib; \etab, \xib) \qquad \forall\,\xib, \etab \in \Hb^1(\Tcal_h)\cap \Lb^\infty.
    \end{equation}
\end{lemma}
\begin{proof}
Using the identity~\eqref{eq:divutv}, the assumption that $\div \psib = 0$, and elementwise integration by parts, we have
\begin{equation}
    \begin{aligned}
        \int_\element (\psib \cdot \gradn) \xib \cdot \etab \dx &= \int_\element \div(\xib\otimes\psib) \cdot \etab - (\div\psib)\xib\cdot\etab \dx\\
        &= -\int_\element \xib \otimes \psib : \grad \etab \dx + \int_{\partial \element} ((\xib \otimes \psib) \nb_\element) \cdot \etab \ds\\
        &= -\int_\element (\psib \cdot \gradn) \etab \cdot \xib \dx + \int_{\partial \element} (\psib \cdot \nb_\element) \xib \cdot \etab \ds,
    \end{aligned}
\end{equation}
where $\nb_\element$ is the unit outward normal vector on $\partial \element$. Summing over all elements $\element \in \Tcal_h$ and noting that $\psib\cdot\nb_\face= 0$ on $\face \in \Fcal^\partial_h$, we have the desired identity
\begin{equation}
    ((\psib \cdot \gradn) \xib, \etab)_{\Tcal_h} = -((\psib \cdot \gradn) \etab, \xib)_{\Tcal_h} + \langle (\psib \cdot \nb) \jump{\xib}, \avg{\etab}\rangle_{\Fcal^\circ_h} + \langle (\psib \cdot \nb) \jump{\etab}, \avg{\xib}\rangle_{\Fcal^\circ_h}.
\end{equation}
\end{proof}
\noindent As is well known, the trilinear form $\tsf_h$ gives rise to a \emph{central scheme}. Suitable stabilization is needed to improve accuracy in the convection-dominated regime when $\nu_\mathrm{v} \ll h$ (see, e.g.,~\cite{brezzi_2004,guzman_2017,barrenechea_2020}). We introduce the standard \emph{interior penalty} stabilization term as
\begin{equation}\label{eq:ssfvel}
    \ssf^\mathrm{v}_h(\xib, \etab) := \gamma_{\mathrm{v}} \langle \jump{\xib}, \jump{\etab}\rangle_{\Fcal^\circ_h}
\end{equation}
where $\gamma_{\mathrm{v}} > 0$ is a penalty parameter. 

\paragraph{Lorentz force and magnetic diffusion.}
Since we use $\Hcurl$-conforming elements for the magnetic induction $\Bb$, the Lorentz force $(\curl \Bb) \times \Bb$ is well-defined; the magnetic diffusion term $\nu_\mathrm{m} \curl \curl \Bb$ also entails no special treatment. 

\paragraph{Magnetic convection.}
The trilinear form $\csf$ in~\eqref{eq:csf} is well-defined on $\Sigmab_h \times \Sigmab_h \times \Vb_{h,0}$. However, like $\tsf_h$ in~\eqref{eq:tsfh}, it gives rise to a central scheme for the magnetic convection term $\curl(\ub \times \Bb)$. Suitable stabilization is needed to improve accuracy in the convection-dominated regime when $\nu_\mathrm{m} \ll h$. We stress that the choice of magnetic stabilization requires particular care, as it plays a central role in our analysis. Specifically, we require it to be
\emph{
\begin{enumerate}[label=$(R\arabic*)$]
    \item strong enough to control the error arising from the magnetic convection term $\curl(\ub \times \Bb)$ \label{eq:req-control-magconv};
    \item weak enough not to hinder the convergence rate;\label{eq:req-weak}
    \item vanishing on $\grad P_h$ to dispense with magnetic Lagrange multipliers.\label{eq:req-weak-divfree}
\end{enumerate}}
\noindent At the moment, the concrete meaning of the above requirements might be fuzzy, but we will revisit them at the concrete steps. We adopt the following curl-jump penalization
\begin{equation}\label{eq:ssfmag}
    \ssf^\mathrm{m}_h(\xib, \etab) := \gamma_{\mathrm{m}} h\langle \jump{\curl\xib}\times \nb, \jump{\curl\etab}\times \nb \rangle_{\Fcal_h}
\end{equation}
where $\gamma_{\mathrm{m}} > 0$ is a penalty parameter. 
\begin{remark}[on the magnetic stabilization $\ssf^\mathrm{m}_h$]\label{rmk:cip}
    The \emph{curl-jump} penalization~\eqref{eq:ssfmag} is a continuous interior penalty (CIP) stabilization, which has been utilized in~\cite[Section~6.2]{daveiga_2026} to stabilize the magnetic convection term. We would like to draw your attention to the scaling of $\ssf^\mathrm{m}_h$. It has a pre-factor of $h$ while that of~\cite[Section~6.2]{daveiga_2026} scales as $h^2$. Since an $h^2$ scaling is typical for CIP-type stabilization (see, e.g.,~\cite{burman_2004}), the choice here leads to an \emph{over-penalization}. The purpose is to provide us with sufficient control over the magnetic convection term~\emph{(Requirement~\ref{eq:req-control-magconv})}. On the other hand, it is affordable thanks to good curl approximation properties of the Whitney $1$-form~\emph{(Requirement~\ref{eq:req-weak})}. 

    Another noteworthy point is that no \emph{value-jump} penalization is imposed. As pointed out in~\cite[Section~6.1]{daveiga_2026}, the value-jump penalization makes the weakly divergence-free space $\Zb^\mathrm{m}_h$ no longer invariant. A magnetic Lagrange multiplier will be needed to enforce the weak divergence-free constraint of $\Bb_h$. In contrast, the curl-jump penalization preserves this constraint over time. Indeed, from the discretization~\eqref{eq:semidiscrete} following, we have $\partial_t(\Bb_h, \grad w_h) = 0$ for all $w_h \in P_h$ since $\ssf^\mathrm{m}_h$ vanishes on $\grad P_h$~\emph{(Requirement~\ref{eq:req-weak-divfree})}.
\end{remark}

\subsection{The semi-discrete scheme}\label{sec:scheme}
We propose the following semi-discrete scheme for the MHD system~\eqref{eq:mhd-fluid}--\eqref{eq:mhd-bc}: 
\begin{formulation}[top=-0.5em]\label{eq:semidiscrete}
    \begin{aligned}
        \text{Seek }(\ub_h(t), p_h(t), \Bb_h(t)) \in \Vb_{h,0} \times S_{h,0} \times \Sigmab_h\text{ such that}  \qquad\qquad\qquad &\\
        (\partial_t \ub_h, \vb_h) + \nu_\mathrm{v} \asf_h(\ub_h, \vb_h) + \tsf_h(\ub_h; \ub_h, \vb_h) + \ssf^\mathrm{v}_h(\ub_h, \vb_h) \qquad\qquad&\\
        - \csf(\Bb_h; \Bb_h, \vb_h) + (p_h, \div\vb_h)  &= (\fb, \vb_h)&&\forall\,\vb_h \in \Vb_{h,0}, \\
        (q_h, \div\ub_h) &= 0 &&\forall\,q_h \in S_{h,0}, \\
        (\partial_t \Bb_h, \Cb_h) + \nu_\mathrm{m} (\curl\Bb_h, \curl\Cb_h) + \csf(\Cb_h; \Bb_h, \ub_h) + \ssf^\mathrm{m}_h(\Bb_h, \Cb_h) &= 0 &&\forall\,\Cb_h \in \Sigmab_h,\\
        \text{with the initial conditions } \ub_h(0) = \IRTavz \ub_0,\; \Bb_h(0) = \Rbcal^\mathrm{ND}_h \Bb_0. \qquad\quad
    \end{aligned}  
\end{formulation}
Here, $\IRTavz$ and $\Rbcal^\mathrm{ND}_h$ are the interpolation operators defined in Section~\ref{sec:cochainproj} and Section~\ref{sec:rh}, respectively. 
For readability, we repeat the definitions of the bilinear and trilinear forms in the scheme~\eqref{eq:semidiscrete} below:
\begin{align}
    \asf_h(\xib, \etab) &:= (\grad \xib, \grad \etab)_{\Tcal_h} - \langle \avg{\grad \xib}\nb, \jump{\etab}\rangle_{\Fcal_h}  \notag\\
    &\qquad\qquad\qquad\qquad\qquad - \langle \avg{\grad \etab}\nb, \jump{\xib}\rangle_{\Fcal_h} + \gamma_\asf h^{-1}\langle \jump{\xib}, \jump{\etab}\rangle_{\Fcal_h},\\
    \tsf_h(\psib; \xib, \etab) &:= ((\psib \cdot \gradn) \xib, \etab)_{\Tcal_h} - \langle (\psib \cdot \nb_\face) \jump{\xib}, \avg{\etab}\rangle_{\Fcal^\circ_h},\\
    \ssf^\mathrm{v}_h(\xib, \etab) &:= \gamma_{\mathrm{v}} \langle \jump{\xib}, \jump{\etab}\rangle_{\Fcal^\circ_h},\\
    \csf(\psib; \xib, \etab) &:= ((\curl \psib) \times \xib, \etab)_{\Tcal_h},\\
    \ssf^\mathrm{m}_h(\xib, \etab) &:= \gamma_{\mathrm{m}} h\langle \jump{\curl\xib}\times \nb, \jump{\curl\etab}\times \nb \rangle_{\Fcal_h}.
\end{align}
\begin{remark}[equivalent formulation]
    Note the vanishing gradient property~\eqref{eq:RTzerograd}. The first equation of the scheme~\eqref{eq:semidiscrete} can be simplified to
    \begin{multline}
        (\partial_t \ub_h, \vb_h) + \nu_\mathrm{v} \gamma_\asf h^{-1}\langle \jump{\ub_h}, \jump{\vb_h}\rangle_{\Fcal_h} - \langle (\ub_h\cdot\nb)\jump{\ub_h},\avg{\vb_h}\rangle_{\Fcal^\circ_h} + \gamma_\mathrm{v}\langle\jump{\ub_h},\jump{\vb_h}\rangle_{\Fcal^\circ_h} \\ - (\curl\Bb_h\times\Bb_h, \vb_h) + (p_h, \div\vb_h)  = (\fb, \vb_h).
    \end{multline}
\end{remark}

\subsection{Divergence constraints, pressure robustness, and energy dissipation}
Since $\div \Vb_{h,0} = S_{h,0}$, the second equation of~\eqref{eq:semidiscrete} implies that the discrete velocity field $\ub_h \in \Zb^{\mathrm{v}}_h$, i.e., it is \emph{strongly} divergence-free. This leads to the \emph{pressure-robustness} (see, e.g.,~\cite{john_2017}) of the scheme~\eqref{eq:semidiscrete} in the sense that any gradient force $\grad \phi$ has no impact on the discrete velocity $\ub_h$. Indeed, one may restrict the velocity trial and test spaces to $\Zb^{\mathrm{v}}_h$ and drop the pressure term. That is, the first equation of~\eqref{eq:semidiscrete} reduces to
\begin{equation}\label{eq:semidiscrete-vel-reduced}
    (\partial_t \ub_h, \vb_h) + \nu_\mathrm{v} \asf_h(\ub_h, \vb_h) + \tsf_h(\ub_h; \ub_h, \vb_h) + \ssf^\mathrm{v}_h(\ub_h, \vb_h)
        - \csf(\Bb_h; \Bb_h, \vb_h) = (\fb, \vb_h)\quad\forall\,\vb_h \in \Zb^{\mathrm{v}}_h.
\end{equation}
The fact that $(\fb + \grad \phi, \vb_h) = (\fb, \vb_h)$ for all $\vb_h \in \Zb^{\mathrm{v}}_h$ implies independence of $\ub_h$ from the external force $\grad \phi$. 

Testing the third equation of~\eqref{eq:semidiscrete} with $\Cb_h = \grad q_h$ for any $q_h \in P_h$ and using the fact that $\curl \circ \grad = \zerob$, we see that $\Bb_h$ is \emph{weakly} divergence-free, i.e., $(\partial_t \Bb_h, \grad q_h) = 0$ for all $q_h \in P_h$ (cf. Requirement~\ref{eq:req-weak-divfree}). Since the interpolation $\Rbcal^\mathrm{ND}_h$ maps functions in $\Zb^{\mathrm{m}}$ to $\Zb^{\mathrm{m}}_h$, the initial condition $\Bb_h(0) = \Rbcal^\mathrm{ND}_h \Bb_0 \in \Zb^{\mathrm{m}}_h$ and thus $\Bb_h(t) \in \Zb^{\mathrm{m}}_h$ for all $t > 0$. Hence, we can restrict both the trial and test spaces of the third equation of~\eqref{eq:semidiscrete} to $\Zb^{\mathrm{m}}_h$ without changing the solution.

Set $\vb_h = \ub_h$ in~\eqref{eq:semidiscrete-vel-reduced} and $\Cb_h = \Bb_h$ in the third equation of~\eqref{eq:semidiscrete}. Sum the two equations. By the skew-symmetry~\eqref{eq:tsfh-skew} of $\tsf_h$ and the cancellation of the Lorentz force term and the magnetic convection term, we have the following energy dissipation law:
\begin{equation}\label{eq:energy-stability}
    \frac{1}{2}\frac{\mathrm{d}}{\mathrm{d}t}\left(\|\ub_h\|^2_{\Lb^2} + \|\Bb_h\|^2_{\Lb^2}\right) + \nu_\mathrm{v} \asf_h(\ub_h, \ub_h) + \nu_\mathrm{m} \|\curl\Bb_h\|^2_{\Lb^2} + \ssf^\mathrm{v}_h(\ub_h, \ub_h) + \ssf^\mathrm{m}_h(\Bb_h, \Bb_h) = (\fb, \ub_h).
\end{equation}

\subsection{Auxiliary trilinear estimates}\label{sec:aux}
In this subsection, we present some auxiliary estimates regarding the trilinear form $\csf$ defined in~\eqref{eq:csf}, which is of central importance in the analysis of the scheme. In particular, we will prove a modified version of the assumption~\cite[Eq.~(6.3)]{daveiga_2026} tailored to our scheme.

The first lemma improves~\cite[Lemma~6.2]{daveiga_2026} by removing the convexity assumption on the domain $\Omega$ using a density argument.
\begin{lemma}[continuous trilinear estimate]\label{lem:conttrilest}
    Given $\ub \in \Hb^1_0 \cap \Wb^{1,\infty}$, there exists a positive constant $C$ such that
    \begin{equation}\label{eq:tril}
        |\csf(\wb; \wb, \ub)| \leq C|\ub|_{\Wb^{1,\infty}}\|\wb\|^2_{\Lb^2}\qquad \forall\,\wb \in \Hcurl\cap\Hdivz.
    \end{equation} 
\end{lemma}
\begin{proof}
    For $\wb \in \Hcurl\cap\Hdivz$, find a sequence $\wb_k \in \Cb^\infty(\overline{\Omega})$ such that $\|\wb_k - \wb\|_{\Hdiv} + \|\wb_k - \wb\|_{\Hcurl} \to 0$ as $k \to \infty$~\cite[Proposition~2.3]{amrouche_1998}. Importantly, $\|\div \wb_k\|_{L^2} \to 0$ as $k \to \infty$. 
    By integration by parts and the vector calculus identity~\eqref{eq:curlutv}, we deduce that
    \begin{equation}
        \begin{aligned}
            \csf(\wb_k; \wb_k, \ub) &:= ((\curl\wb_k) \times \wb_k, \ub) = (\wb_k \times \ub, \curl\wb_k) = (\curl(\wb_k \times \ub), \wb_k) \\
            &= ((\div \ub)\wb_k, \wb_k) - ((\div \wb_k)\ub, \wb_k) + ((\ub \cdot \gradn)\wb_k, \wb_k) - ((\wb_k \cdot \gradn)\ub, \wb_k). 
        \end{aligned}
    \end{equation}
    The first and last terms are bounded by $C|\ub|_{\Wb^{1,\infty}}\|\wb_k\|^2_{\Lb^2}$. To estimate the third term, we conduct integration by parts again
    \begin{equation}\label{eq:ugradwkwk}
        ((\ub \cdot \gradn)\wb_k, \wb_k) = \frac{1}{2}(\ub, \grad(\wb_k \cdot \wb_k)) = -\frac{1}{2}((\div \ub)\wb_k, \wb_k) + \frac{1}{2}\langle (\ub\cdot\nb)\wb_k, \wb_k\rangle_{\Gamma}.
    \end{equation}
    The boundary term vanishes since $\ub \in \Hb^1_0$ and the left-hand side of~\eqref{eq:ugradwkwk} is bounded by $C|\ub|_{\Wb^{1,\infty}}\|\wb_k\|^2_{\Lb^2}$. So far, we have shown that 
    \begin{equation}
        |\csf(\wb_k; \wb_k, \ub)| \leq C|\ub|_{\Wb^{1,\infty}}\|\wb_k\|^2_{\Lb^2} + |((\div \wb_k)\ub, \wb_k)|.
    \end{equation}
    Pass to the limit $k \to \infty$, and the desired result follows.
\end{proof}

To motivate the discrete estimate, recall the assumption~\cite[Eq.~(6.3)]{daveiga_2026} stating that 
\begin{equation}\label{eq:daveiga-assumption}
    \csf(\wb_h; \wb_h, \ub) \leq C|\ub|_{\Wb^{1,\infty}}\|\wb_h\|^2_{\Lb^2} \qquad \forall\,\ub \in \Hb^1_0 \cap \Wb^{1,\infty},\;\wb_h \in \Zb^{\mathrm{m}}_h.
\end{equation}
This can be interpreted as a discrete version of Lemma~\ref{lem:conttrilest}.
Since $\wb_h \in \Zb^{\mathrm{m}}_h$ is not necessarily divergence-free, the proof of Lemma~\ref{lem:conttrilest} does not apply directly. The idea is to connect the discrete function $\wb_h \in \Zb^{\mathrm{m}}_h$ with a function $\wb \in \Hcurl\cap\Hdivz$ and then apply the continuous trilinear estimate~\eqref{eq:tril}.
Such a connection is provided by the lifting map $\Ksf$ introduced in the next lemma.
\begin{lemma}[a lifting map]\label{lem:ksf}
    There exist a linear map $\Ksf: \Lb^2 \to H^1$ and a positive constant $C$ independent of $h$ such that 
    \begin{equation}\label{eq:ksf-div0}
        \wb - \grad \Ksf\wb \in \Hdivz \qquad \forall\,\wb \in \Lb^2,
    \end{equation}
    and 
    \begin{align}
        \qquad\qquad\qquad\qquad\|\grad \Ksf\wb\|_{L^2} &\leq C\|\wb\|_{\Lb^2}&&\forall\,\wb \in \Lb^2,\qquad\qquad \label{eq:ksf-stabh1}\\
        \qquad\qquad\qquad\qquad\|\Ksf\wb\|_{L^2} &\leq C h\|\wb\|_{\Lb^2}&&\forall\,\wb \in \Lb^2 \cap (\grad P_h)^\perp.\qquad\qquad\label{eq:ksf-stabl2}
    \end{align}
\end{lemma}
\begin{proof}
Let $\Omega^e \supseteq \Omega$ be a \emph{convex} polyhedral domain containing the original polyhedron. Given $\wb \in \Lb^2$, we denote its zero extension to $\Omega^e$ as $\wb^e$. Then, consider the following Dirichlet problem on the extended domain $\Omega^e$: Seek $\varphi_\wb \in H^1_0(\Omega^e)$ such that
\begin{equation}\label{eq:varphipoisson-extended}
    (\grad \varphi_\wb, \grad q)_{\Omega^e} = (\wb^e, \grad q)_{\Omega^e} \qquad \forall\,q \in H^1_0(\Omega^e).
\end{equation}
The Lax--Milgram theorem gives $\|\varphi_\wb\|_{H^1(\Omega^e)} \leq C\|\wb^e\|_{\Lb^2(\Omega^e)} = C\|\wb\|_{\Lb^2}$. 
Thereby, we define the map 
\begin{equation}\label{eq:ksfdef}
    \Ksf: \Lb^2 \ni \wb \mapsto \varphi_\wb|_\Omega \in H^1
\end{equation}
where $\varphi_\wb$ is the solution of~\eqref{eq:varphipoisson-extended}. The estimate~\eqref{eq:ksf-stabh1} is obvious. Since~\eqref{eq:varphipoisson-extended} implies that $\div(\wb^e - \grad \varphi_\wb) = 0$ in $\Omega^e$ distributionally, we have $\div(\wb - \grad \Ksf\wb) = 0$ in $\Omega$ distributionally, which verifies~\eqref{eq:ksf-div0}. To show~\eqref{eq:ksf-stabl2}, we apply the standard duality argument. Let $\eta \in H^1_0(\Omega^e)$ solve the dual problem
\begin{equation}\label{eq:dualpoisson-extended}
    (\grad \eta, \grad q)_{\Omega^e} = (\varphi_\wb, q)_{\Omega^e} \qquad \forall\,q \in H^1_0(\Omega^e).
\end{equation}
Since $\Omega^e$ is convex, by elliptic regularity, we have $\eta \in H^2(\Omega^e)$ and $\|\eta\|_{H^2(\Omega^e)} \leq C\|\varphi_\wb\|_{L^2(\Omega^e)}$. Set $q = \varphi_\wb$ in~\eqref{eq:dualpoisson-extended} and let $\eta_I \in P_h$ be the nodal interpolant of $\eta|_\Omega$. We have
\begin{equation}
    \begin{aligned}
        \|\varphi_\wb\|^2_{L^2(\Omega^e)} &= (\grad \eta, \grad \varphi_\wb)_{\Omega^e} = (\grad \eta, \wb^e)_{\Omega^e} = (\grad \eta, \wb) \\
        &=(\grad (\eta - \eta_I), \wb) \leq Ch\|\eta\|_{H^2(\Omega)}\|\wb\|_{\Lb^2} \leq Ch\|\varphi_\wb\|_{L^2(\Omega^e)}\|\wb\|_{\Lb^2},
    \end{aligned}
\end{equation}
where we used the condition that $\wb \in (\grad P_h)^\perp$ to get $(\grad \eta_I, \wb)_{\Omega} = 0$. We can conclude~\eqref{eq:ksf-stabl2}.
\end{proof}
\noindent We remark that the enlarged convex domain $\Omega^e$ provides the $H^2$-regularity needed in the above duality argument without assuming $\Omega$ to be convex.

Finally, we are ready to show a discrete counterpart of Lemma~\ref{lem:conttrilest} for $\wb_h \in \Zb^{\mathrm{m}}_h$ (cf.~\eqref{eq:daveiga-assumption}).
\begin{lemma}[discrete trilinear estimate]\label{lem:disctrilest}
    Given $\ub \in \Hb^1_0 \cap \Wb^{1,\infty}$, there exists a positive constant $C$ independent of $h$ such that
    \begin{equation}\label{eq:disctrilest}
        |\csf(\wb_h; \wb_h, \ub)| \leq C|\ub|_{\Wb^{1,\infty}}\|\wb_h\|_{\Lb^2}^2 + Ch^{1/2}\|\ub\|_{\Lb^{\infty}}\|\wb_h\|_{\Lb^2}\left\|\jump{\curl\wb_h}\times\nb\right\|_{\Fcal^\circ_h} \quad\forall\,\wb_h \in \Zb^{\mathrm{m}}_h.
    \end{equation}
\end{lemma}
\begin{proof}
    Let $\wb := \wb_h - \grad \Ksf\wb_h$ where $\Ksf$ is the potential mapping defined in Lemma~\ref{lem:ksf}. It holds that $\wb \in \Hcurl\cap\Hdivz$. 
    Then, we deduce
    \begin{equation}
        \csf(\wb_h; \wb_h, \ub) := (\curl\wb_h \times \wb_h, \ub) = (\curl\wb \times \wb, \ub) + (\curl\wb_h \times \grad \Ksf\wb_h, \ub).
    \end{equation}
    The first term is bounded by $C|\ub|_{\Wb^{1,\infty}}\|\wb_h\|^2_{\Lb^2}$ due to Lemma~\ref{lem:conttrilest}. To estimate the second term, we perform elementwise integration by parts and deduce that
    \begin{equation}\label{eq:curlwgradKsf}
        \begin{aligned}
            (\curl\wb_h \times \grad \Ksf\wb_h, \ub) &= (\ub\times \curl \wb_h, \grad \Ksf\wb_h) \\
            &= - (\div(\ub\times\curl\wb_h), \Ksf\wb_h)_{\Tcal_h} + \langle \ub\times\jump{\curl\wb_h}\cdot\nb, \avg{\Ksf\wb_h}\rangle_{\Fcal^\circ_h}.
        \end{aligned}
    \end{equation}
    Here, we used the fact that the interior jump $\jump{\Ksf\wb_h}_\face = 0$ for all $\face \in \Fcal^\circ_h$ and $\ub|_\Gamma = \zerob$. Using the identity~\eqref{eq:divutimesv} and noticing that $\curl \curl \wb_h = \zerob$ elementwise, the first term on the right-hand side of~\eqref{eq:curlwgradKsf} is bounded by 
    \begin{equation}
        \left|(\div(\ub\times\curl\wb_h), \Ksf\wb_h)_{\Tcal_h}\right| = \left|((\curl \ub)\cdot \curl \wb_h, \Ksf\wb_h)_{\Tcal_h}\right| \leq C|\ub|_{\Wb^{1,\infty}}\|\wb_h\|^2_{\Lb^2}
    \end{equation}
    where we used the inverse inequality and the stability estimate~\eqref{eq:ksf-stabl2}. To estimate the second term on the right-hand side of~\eqref{eq:curlwgradKsf}, we invoke the multiplicative trace inequality (see, e.g.,~\cite[Eq.~(12.17)]{ern_2021}):
    \begin{equation}\label{eq:multiplicativetrace}
        \|\vb\|_{\Lb^2(\face)} \leq C\left(h^{-1/2}\|\vb\|_{\Lb^2(\element)} + h^{s-1/2}|\vb|_{\Hb^s(\element)}\right)\qquad \forall \,\vb \in \Hb^s(\element),\;\face \subset \partial \element,\;s \in (1/2, 1].
    \end{equation}
    Then, we deduce that
    \begin{equation}
        \begin{aligned}
        \left|\langle \ub\times \jump{\curl\wb_h}\cdot\nb, \avg{\Ksf\wb_h}\rangle_{\Fcal^\circ_h}\right| &= \left|\langle \jump{\curl\wb_h}\times\nb\cdot \ub, \avg{\Ksf\wb_h}\rangle_{\Fcal^\circ_h}\right| \\
        &\leq C\|\ub\|_{\Lb^\infty}\left\|\jump{\curl\wb_h}\times \nb \right\|_{\Fcal^\circ_h}\left(h^{-1/2}\|\Ksf\wb_h\|_{L^2} + h^{1/2}|\Ksf\wb_h|_{H^1}\right) \\
        &\leq Ch^{1/2}\|\ub\|_{\Lb^\infty}\left\|\jump{\curl\wb_h}\times \nb \right\|_{\Fcal^\circ_h}\|\wb_h\|_{\Lb^2}
        \end{aligned}
    \end{equation}
    where we used the stability estimates~\eqref{eq:ksf-stabh1} and~\eqref{eq:ksf-stabl2}.
    The desired result follows from the above estimates.
\end{proof}
\begin{remark}[comparison of~\eqref{eq:daveiga-assumption} and~\eqref{eq:disctrilest}]
    In contrast to the original assumption~\eqref{eq:daveiga-assumption}, the new estimate~\eqref{eq:disctrilest} is weaker, involving an additional jump term $h^{1/2}\|\jump{\curl\wb_h}\times\nb\|_{\Fcal^\circ_h}$ on the right-hand side. 
\end{remark}

The last lemma in this subsection is a technical lemma stating a variant of Lemma~\ref{lem:disctrilest}. Since the proof contains no new ideas beyond those of Lemmas~\ref{lem:conttrilest} and~\ref{lem:disctrilest}, we merely sketch it.
\begin{lemma}[discrete trilinear estimate II]\label{lem:disctrilest2}
Given $\ub \in \Hb^1_0 \cap \Wb^{1,\infty}$, there exists a positive constant $C$ independent of $h$ such that
\begin{equation}\label{eq:disctrilest2}
    |\csf(\wb_h; \vb, \ub)| \leq C\|\ub\|_{\Wb^{1,\infty}}\|\wb_h\|_{\Lb^2}\|\vb\|_{\Hcurl} + Ch^{1/2}\|\ub\|_{\Lb^{\infty}}\|\vb\|_{\Lb^2}\left\|\jump{\curl\wb_h}\times\nb\right\|_{\Fcal^\circ_h}
\end{equation}
for all $\wb_h \in \Zb^{\mathrm{m}}_h$ and $\vb \in \Hcurl\cap (\grad P_h)^\perp$.
\end{lemma} 
\begin{proof}[Sketch of proof]
    Consider the decomposition $\wb_h = (\wb_h - \grad \Ksf\wb_h) + \grad \Ksf\wb_h$ and $\vb = (\vb - \grad \Ksf\vb) + \grad \Ksf\vb$. The first components of both terms are amenable to the arguments in Lemma~\ref{lem:conttrilest}. Since both $\wb_h$ and $\vb$ are orthogonal to $\grad P_h$, the second components can be effectively controlled by the stability estimates~\eqref{eq:ksf-stabh1} and~\eqref{eq:ksf-stabl2}. Repeat the derivation in Lemma~\ref{lem:disctrilest} and the desired result follows.
\end{proof}

\begin{remark}[rationale of the trilinear estimates]
    Conceptually, the trilinear estimates in Lemmas~\ref{lem:disctrilest} and~\ref{lem:disctrilest2} transfer the curl operator away from the discrete function $\wb_h$ at the cost of an additional curl jump term, which motivates~$\ssf^\mathrm{m}_h$~\eqref{eq:ssfmag} in our scheme~\eqref{eq:semidiscrete}. 
    In this way, one avoids the $O(h^{-1})$-factor from applying the inverse inequality, which is crucial for an improved convergence rate in the convection-dominated regime $\nu_\mathrm{m} \ll h$. 
\end{remark}

\subsection{A priori error estimate}\label{sec:err-est}
Now we present the main result on a priori error estimates of the semi-discrete scheme~\eqref{eq:semidiscrete}.
\begin{theorem}[a priori error estimate]\label{lem:err-est}
    Let $(\ub, \Bb)$ solve the MHD system~\eqref{eq:MHD-weak} and $(\ub_h, \Bb_h) \in \Vb_{h,0} \times \Sigmab_h$ solve the semi-discrete scheme~\eqref{eq:semidiscrete}. Suppose Assumption~\ref{ass:reg} holds. Then, there exists a positive constant $C$ depending on the solution norms (specified in Assumption~\ref{ass:reg}) and the penalty parameters $\gamma_\asf, \gamma_\mathrm{v}, \gamma_\mathrm{m},$ but independent of $h$, $\nu_\mathrm{v}$, and $\nu_\mathrm{m}$ such that
    \begin{equation}\label{eq:err-est}
        \sup_{t \in (0,T)}\left(\|\ub(t) - \ub_h(t)\|_{\Lb^2} + \|\Bb(t) - \Bb_h(t)\|_{\Lb^2}\right) \leq C (h^{1/2} + \nu^{1/2}_\mathrm{v} + \nu^{1/2}_\mathrm{m} h^{1/2}).
    \end{equation}
\end{theorem}

\begin{proof}
First, from the derivation in Section~\ref{sec:indivterms}, we see that the scheme~\eqref{eq:semidiscrete} is consistent with the exact solution $(\ub, \Bb)$, that is,
\begin{equation}\label{eq:consistency}
\begin{aligned}
    (\partial_t \ub, \vb_h) + \nu_\mathrm{v} \asf_h(\ub, \vb_h) + \tsf_h(\ub; \ub, \vb_h) + \ssf^\mathrm{v}_h(\ub, \vb_h) - \csf(\Bb; \Bb, \vb_h) &= (\fb, \vb_h) &&\forall\,\vb_h \in \Zb^{\mathrm{v}}_h,\\
    (\partial_t \Bb, \Cb_h) + \nu_\mathrm{m} (\curl\Bb, \curl\Cb_h) + \csf(\Cb_h; \Bb, \ub) + \ssf^\mathrm{m}_h(\Bb, \Cb_h) &= 0 &&\forall\,\Cb_h \in \Zb^{\mathrm{m}}_h.
\end{aligned}
\end{equation}
Indeed, under the regularity stated in Assumption~\ref{ass:reg}, it is clear that all the terms in~\eqref{eq:consistency} are well-defined. By definition of $\tsf_h, \ssf^\mathrm{v}_h, \ssf^\mathrm{m}_h$ in~\eqref{eq:tsfh},~\eqref{eq:ssfvel}, and~\eqref{eq:ssfmag}, it holds that  $\tsf_h(\ub; \ub, \vb_h) = \tsf(\ub; \ub, \vb_h)$, $\ssf_h^\mathrm{v}(\ub, \vb_h) = 0$, and $\ssf_h^\mathrm{m}(\Bb, \Cb_h) = 0$.
In view of the continuous formulation~\eqref{eq:MHD-weak} and the consistency~\eqref{eq:vecdiffusion-consistency} of $\asf_h$, we arrive at~\eqref{eq:consistency}.

Denote the interpolants $\ub_I := \IRTavz \ub,\;\Bb_I := \Rbcal^\mathrm{ND}_h \Bb$ and the error terms
\begin{equation}
    \eb^\mathrm{v}_h := \ub_h - \ub_I, \qquad \eb^\mathrm{v}_I := \ub - \ub_I, \qquad \eb^\mathrm{m}_h := \Bb_h - \Bb_I,\qquad \eb^\mathrm{m}_I := \Bb - \Bb_I.
\end{equation}
Combining the scheme~\eqref{eq:semidiscrete} with the consistency~\eqref{eq:consistency}, we obtain the following ``Galerkin orthogonality'':
\begin{equation}\label{eq:galorth}
    \begin{aligned}
        (\partial_t (\ub - \ub_h), \vb_h) + \nu_\mathrm{v} \asf_h(\ub - \ub_h, \vb_h)  + \ssf^\mathrm{v}_h(\ub - \ub_h, \vb_h)\qquad\qquad\qquad\qquad&\\
        + \tsf_h(\ub; \ub, \vb_h) - \tsf_h(\ub_h; \ub_h, \vb_h)- \csf(\Bb; \Bb, \vb_h) + \csf(\Bb_h; \Bb_h, \vb_h) &= 0, \\
        (\partial_t (\Bb - \Bb_h), \Cb_h) + \nu_\mathrm{m} (\curl(\Bb - \Bb_h), \curl\Cb_h) + \ssf^\mathrm{m}_h(\Bb- \Bb_h, \Cb_h)  \qquad\qquad&\\
        + \csf(\Cb_h; \Bb, \ub) - \csf(\Cb_h; \Bb_h, \ub_h) &= 0.
    \end{aligned}
\end{equation}
Setting $\vb_h = \eb^\mathrm{v}_h \in \Zb^{\mathrm{v}}_h$ and $\Cb_h = \eb^\mathrm{m}_h \in \Zb^{\mathrm{m}}_h$ in~\eqref{eq:galorth}, we get
\begin{multline}\label{eq:err-eq}
    \frac{1}{2}\frac{\mathrm{d}}{\mathrm{d}t}\left(\|\eb^\mathrm{v}_h\|^2_{\Lb^2} + \|\eb^\mathrm{m}_h\|^2_{\Lb^2}\right) + \nu_\mathrm{v} \asf_h(\eb^\mathrm{v}_h, \eb^\mathrm{v}_h) + \nu_\mathrm{m} \|\curl\eb^\mathrm{m}_h\|^2_{\Lb^2} \\+ \ssf^\mathrm{v}_h(\eb^\mathrm{v}_h, \eb^\mathrm{v}_h) + \ssf^\mathrm{m}_h(\eb^\mathrm{m}_h, \eb^\mathrm{m}_h) = 
    I^{\mathrm{bil}} + I^{\mathrm{vc}} + I^{\mathrm{lor}} + I^{\mathrm{mc}}
\end{multline}
where the bilinear term $I^{\mathrm{bil}}$, the velocity convection trilinear term $I^{\mathrm{vc}}$, the Lorentz force trilinear term $I^{\mathrm{lor}}$, and the magnetic convection trilinear term $I^{\mathrm{mc}}$ are defined as
\begin{align}
        I^{\mathrm{bil}} &:= (\partial_t \eb^\mathrm{v}_I, \eb^\mathrm{v}_h) + (\partial_t \eb^\mathrm{m}_I, \eb^\mathrm{m}_h) + \nu_\mathrm{v} \asf_h(\eb^\mathrm{v}_I, \eb^\mathrm{v}_h) + \nu_\mathrm{m} (\curl\eb^\mathrm{m}_I, \curl\eb^\mathrm{m}_h) + \ssf^\mathrm{v}_h(\eb^\mathrm{v}_I, \eb^\mathrm{v}_h) + \ssf^\mathrm{m}_h(\eb_I^{\mathrm{m}}, \eb^\mathrm{m}_h),\label{eq:Ibil}\\
        I^{\mathrm{vc}} &:= \tsf_h(\ub; \ub, \eb^\mathrm{v}_h) - \tsf_h(\ub_h; \ub_h, \eb^\mathrm{v}_h),\label{eq:Ivc}\\
        I^{\mathrm{lor}} &:= \csf(\Bb_I; \Bb_h, \eb^\mathrm{v}_h)-\csf(\Bb; \Bb, \eb^\mathrm{v}_h) ,\label{eq:Ilor}\\
        I^{\mathrm{mc}} &:= \csf(\eb^\mathrm{m}_h; \Bb, \ub) - \csf(\eb^\mathrm{m}_h; \Bb_h, \ub_I).\label{eq:Imc}
\end{align}
Here, we stress that the derivation of~\eqref{eq:Ilor},~\eqref{eq:Imc} involves the following cancellation identity (cf.~the derivation of the energy stability~\eqref{eq:energy-stability}):
\begin{equation}\label{eq:cancellation}
    \csf(\Bb_h; \Bb_h, \eb^\mathrm{v}_h) - \csf(\eb^\mathrm{m}_h; \Bb_h, \ub_h) = \csf(\Bb_I; \Bb_h, \eb^\mathrm{v}_h) - \csf(\eb^\mathrm{m}_h; \Bb_h, \ub_I).
\end{equation}

\begin{remark}[general strategy behind the error analysis]
    One should avoid trilinear terms whose three arguments all involve the discrete solutions, i.e., $\ub_h$ and $\Bb_h$, which are in general hard to control. This is where the cancellation~\eqref{eq:cancellation} comes into play. The trilinear terms with discrete solutions in two arguments are to be absorbed into the left-hand side of~\eqref{eq:err-eq}. The remaining terms are controlled by interpolation approximations.
\end{remark}

In what follows, we treat the four terms $I^{\mathrm{bil}}$, $I^{\mathrm{vc}}$, $I^{\mathrm{lor}}$, and $I^{\mathrm{mc}}$ one by one. 

\paragraph{Bilinear term $I^{\mathrm{bil}}$.}
All the terms in $I^{\mathrm{bil}}$ that do not involve face integrals are directly estimated by the Cauchy--Schwarz inequality and the interpolation approximations in Sections~\ref{sec:cochainproj} and~\ref{sec:rh}. We have
\begin{align}
    (\partial_t \eb^\mathrm{v}_I, \eb^\mathrm{v}_h) &\leq Ch^{1/2}\|\partial_t \ub\|_{\Hb^{1/2+\delta}}\|\eb^\mathrm{v}_h\|_{\Lb^2},\\
    (\partial_t \eb^\mathrm{m}_I, \eb^\mathrm{m}_h) &\leq Ch^{1/2}\|\partial_t \Bb\|_{\Hb^{1/2}}\|\eb^\mathrm{m}_h\|_{\Lb^2},\\
    \nu_\mathrm{m}(\curl\eb^\mathrm{m}_I, \curl\eb^\mathrm{m}_h) &\leq C\nu_\mathrm{m} h^{1/2}|\curl\Bb|_{\Hb^{1/2}}\|\curl\eb^\mathrm{m}_h\|_{\Lb^2}.
\end{align}
To estimate the terms $\nu_\mathrm{v} \asf_h(\eb^\mathrm{v}_I, \eb^\mathrm{v}_h)$ and $\ssf^\mathrm{v}_h(\eb^\mathrm{v}_I, \eb^\mathrm{v}_h)$, we invoke the multiplicative trace inequality~\eqref{eq:multiplicativetrace} and deduce that
\begin{align}
    \nu_\mathrm{v} \asf_h(\eb^\mathrm{v}_I, \eb^\mathrm{v}_h) &\leq C\nu_\mathrm{v} (\|\ub\|_{\Hb^{1}} + h^{1/2+\delta}\|\ub\|_{\Hb^{3/2+\delta}})(\asf_h(\eb^\mathrm{v}_h,\eb^\mathrm{v}_h))^{1/2},\label{eq:nuvasfebviebvh}\\
    \ssf^\mathrm{v}_h(\eb^\mathrm{v}_I, \eb^\mathrm{v}_h) &\leq C h^{1/2}\|\ub\|_{\Hb^1}(\ssf^\mathrm{v}_h(\eb^\mathrm{v}_h, \eb^\mathrm{v}_h))^{1/2},\\
    \ssf^\mathrm{m}_h(\eb_I^{\mathrm{m}}, \eb^\mathrm{m}_h) &\leq C h^{1/2}\left(\|\curl\Bb\|_{\Hb^{1/2+\delta}} + \|\curl\Bb\|_{\Lb^\infty}\right) (\ssf^\mathrm{m}_h(\eb^\mathrm{m}_h, \eb^\mathrm{m}_h))^{1/2}.\label{eq:ssfeiehm2}
\end{align}
Here, we used Lemma~\ref{lem:asf} to obtain~\eqref{eq:nuvasfebviebvh}, and treated $\face \in \Fcal^\partial_h$ separately using the boundary condition $\curl\Bb\times\nb = \zerob$ in~\eqref{eq:ssfeiehm2}.
Revisit Requirement~\ref{eq:req-weak} which states that the magnetic stabilization $\ssf^\mathrm{m}_h$ should not be so strong as to hinder convergence. This is reflected in the above estimate~\eqref{eq:ssfeiehm2}.

Collecting the above estimates, we have
\begin{multline}\label{eq:Ibilest}
    I^{\mathrm{bil}} \leq C_\mathrm{bil} \left(h + \nu_\mathrm{v} + \nu_\mathrm{m} h\right) + \|\eb^\mathrm{v}_h\|^2_{\Lb^2} + \|\eb^\mathrm{m}_h\|^2_{\Lb^2} +  \frac{\nu_\mathrm{v}}{2}\asf_h(\eb^\mathrm{v}_h, \eb^\mathrm{v}_h) + \frac{\nu_\mathrm{m}}{2}\|\curl\eb^\mathrm{m}_h\|^2_{\Lb^2} \\+ \frac{1}{4}\ssf^\mathrm{v}_h(\eb^\mathrm{v}_h, \eb^\mathrm{v}_h) + \frac{1}{4}\ssf^\mathrm{m}_h(\eb^\mathrm{m}_h, \eb^\mathrm{m}_h)
\end{multline}
where $C_\mathrm{bil}(t)$ depends on $\|\partial_t \ub\|_{\Hb^{1/2+\delta}}$, $\|\partial_t \Bb\|_{\Hb^{1/2}}$, $\|\ub\|_{\Hb^{3/2+\delta}}$, $\|\curl\Bb\|_{\Hb^{1/2+\delta}}$, and $\|\curl\Bb\|_{\Lb^\infty}$.

\paragraph{Velocity convection term $I^{\mathrm{vc}}$.}

In the context of $\Hdiv$-conforming discretization of the Navier--Stokes equations, a key step towards an $O(h^{1/2})$-improved $\Lb^2$-error estimate is to invoke the interior moment orthogonality~\cite[Eqs.~(7.8) and~(7.9)]{arnold_2018} of the canonical interpolation; see, e.g.,~\cite[Eq.~(4.15)]{barrenechea_2020} and~\cite[Term~$T_{2,2}$]{daveiga_2025}. However, this approach does not work for the Whitney form $\Vb_h$ due to the absence of the interior degrees of freedom. Instead, we utilize the vanishing-gradient property~\eqref{eq:RTzerograd} to drop the volume integral term and only deal with the face integral term. 

We split $I^{\mathrm{vc}}$ as
\begin{equation}\label{eq:Ivc-split}
    I^{\mathrm{vc}} = \underbrace{\tsf_h(\eb^\mathrm{v}_I; \ub, \eb^\mathrm{v}_h)}_{I^{\mathrm{vc}}_1} - \underbrace{\tsf_h(\ub_I; \eb^\mathrm{v}_h, \eb^\mathrm{v}_I)}_{I^{\mathrm{vc}}_2} - \underbrace{\tsf_h(\eb^\mathrm{v}_h; \ub_I, \eb^\mathrm{v}_h)}_{I^{\mathrm{vc}}_3},
\end{equation}
where we used the skew-symmetry~\eqref{eq:tsfh-skew} multiple times. Since $\ub \in \Hb^1$, the face integrals are zero in $I^{\mathrm{vc}}_1$, whose volume integral is then directly bounded as
\begin{equation}
    I^{\mathrm{vc}}_1 \leq \|\eb^\mathrm{v}_I\|_{\Lb^2} |\ub|_{\Wb^{1,\infty}}\|\eb^\mathrm{v}_h\|_{\Lb^2} \leq C h|\ub|_{\Hb^{1}}|\ub|_{\Wb^{1,\infty}}\|\eb^\mathrm{v}_h\|_{\Lb^2}
\end{equation}
by the $\Lb^2$-approximation~\eqref{eq:approxl2RT0}. The term $I^{\mathrm{vc}}_2$, as we discussed above, is usually amenable to the moment orthogonality of the canonical interpolation for higher-order elements. Yet here, both $I^{\mathrm{vc}}_2$ and $I^{\mathrm{vc}}_3$ are left with face integrals due to the fact that $\eb^\mathrm{v}_h, \ub_I \in \Zb^\mathrm{v}_h$ which therefore satisfy the vanishing-gradient property~\eqref{eq:RTzerograd}. Thereby, we deduce
\begin{equation}
    |I^{\mathrm{vc}}_2| = |\langle (\ub_I \cdot \nb) \jump{\eb^\mathrm{v}_h}, \avg{\eb^\mathrm{v}_I}\rangle_{\Fcal^\circ_h}| \leq Ch^{1/2} \|\ub\|_{\Lb^\infty}|\ub|_{\Hb^1}\|\jump{\eb^\mathrm{v}_h}\|_{\Lb^2(\Fcal^\circ_h)}
\end{equation}
where we used the multiplicative trace inequality~\eqref{eq:multiplicativetrace}, the $\Lb^\infty$-stability from~\eqref{eq:approxlinfRT0}, and the $\Lb^2$-approximation~\eqref{eq:approxl2RT0}. The term $I^{\mathrm{vc}}_3$ is a nonlinear error term and can be bounded as
\begin{equation}
    |I^{\mathrm{vc}}_3| = |\langle (\eb^\mathrm{v}_h \cdot \nb) \jump{\ub_I}, \avg{\eb^\mathrm{v}_h}\rangle_{\Fcal^\circ_h}| = |\langle (\eb^\mathrm{v}_h \cdot \nb) \jump{\ub_I - \ub}, \avg{\eb^\mathrm{v}_h}\rangle_{\Fcal^\circ_h}| \leq C |\ub|_{\Wb^{1,\infty}} \|\eb^\mathrm{v}_h\|^2_{\Lb^2}
\end{equation}
where we used the $\Lb^\infty$-approximation~\eqref{eq:approxlinfRT0}. Combining the above estimates, using Young's inequality, and recalling the definition of $\ssf^\mathrm{v}_h$, we have
\begin{equation}\label{eq:Ivcest}
    I^{\mathrm{vc}} \leq C_\mathrm{vc}(t) h + C_\mathrm{vc}^{\ast}(t)\|\eb^\mathrm{v}_h\|^2_{\Lb^2} + \frac{1}{4}\ssf^\mathrm{v}_h(\eb^\mathrm{v}_h, \eb^\mathrm{v}_h)
\end{equation}
where $C_\mathrm{vc}(t)$ and  $C_\mathrm{vc}^{\ast}(t)$ depend on $\|\ub\|_{\Wb^{1,\infty}}$.

\paragraph{Lorentz force term $I^{\mathrm{lor}}$.}
Estimating $I^{\mathrm{lor}}$ is relatively straightforward. It splits into
\begin{equation}
    I^{\mathrm{lor}} = \underbrace{\csf(\Bb_I; \eb^\mathrm{m}_h, \eb^\mathrm{v}_h)}_{I^{\mathrm{lor}}_1} - \underbrace{\csf(\eb^\mathrm{m}_I; \Bb, \eb^\mathrm{v}_h)}_{I^{\mathrm{lor}}_2} - \underbrace{\csf(\Bb_I; \eb^\mathrm{m}_I, \eb^\mathrm{v}_h)}_{I^{\mathrm{lor}}_3}.
\end{equation}
The $\Lb^\infty$-stability from~\eqref{eq:approxcurlinfR} gives 
\begin{equation}
    I^{\mathrm{lor}}_1 \leq C_1 \|\curl\Bb\|_{\Lb^{\infty}} \|\eb^\mathrm{m}_h\|_{\Lb^2}\|\eb^\mathrm{v}_h\|_{\Lb^2}.
\end{equation}
Both $I^{\mathrm{lor}}_2$ and $I^{\mathrm{lor}}_3$ are bounded by interpolation errors. The $\Hcurl$-approximation~\eqref{eq:approxcurlR} yields
\begin{equation}
    |I^{\mathrm{lor}}_2| \leq C_2 h^{1/2} |\curl \Bb|_{\Hb^{1/2}} \|\Bb\|_{\Lb^\infty} \|\eb^\mathrm{v}_h\|_{\Lb^2},
\end{equation}
and the $\Lb^\infty$-stability from~\eqref{eq:approxcurlinfR} and the $\Lb^2$-approximation~\eqref{eq:approxl2R} give
\begin{equation}
    |I^{\mathrm{lor}}_3| \leq C_3 h^{1/2} \|\curl\Bb\|_{\Lb^\infty} |\Bb|_{\Hb^{1/2}} \|\eb^\mathrm{v}_h\|_{\Lb^2}.
\end{equation}
Combining the above estimates and using Young's inequality, we have
\begin{equation}\label{eq:Ilorest}
    I^{\mathrm{lor}} \leq C_\mathrm{lor} h + C_\mathrm{lor}^\ast\|\eb^\mathrm{m}_h\|^2_{\Lb^2} + C_\mathrm{lor}^\dagger \|\eb^\mathrm{v}_h\|^2_{\Lb^2}
\end{equation}
where $C_\mathrm{lor}(t)$, $C_\mathrm{lor}^\ast(t)$, and $C_\mathrm{lor}^\dagger(t)$ depend on $\|\Bb\|_{\Lb^\infty}$, $\|\Bb\|_{\Hb^{1/2}}$, $\|\curl\Bb\|_{\Lb^\infty}$, and $|\curl\Bb|_{\Hb^{1/2}}$.

\paragraph{Magnetic convection term $I^{\mathrm{mc}}$.} The main burden of the analysis is to estimate the term $I^{\mathrm{mc}}$. We split it into
\begin{equation}
    I^{\mathrm{mc}} = \underbrace{\csf(\eb^\mathrm{m}_h; \Bb, \eb^\mathrm{v}_I)}_{I^{\mathrm{mc}}_1} + \underbrace{\csf(\eb^\mathrm{m}_h; \eb^\mathrm{m}_I, \ub_I)}_{I^{\mathrm{mc}}_2} - \underbrace{\csf(\eb^\mathrm{m}_h; \eb^\mathrm{m}_h, \ub_I)}_{I^{\mathrm{mc}}_3}.
\end{equation}
The term $I^{\mathrm{mc}}_1$ has a similar convection nature to $I^{\mathrm{vc}}_2$ in~\eqref{eq:Ivc-split} and the stabilization term $\ssf^\mathrm{m}_h$ comes into play (see Requirement~\ref{eq:req-control-magconv}). Applying the curl inverse inequality~\eqref{eq:curlinv}, we obtain
\begin{equation}\label{eq:Imcest1}
    I^{\mathrm{mc}}_1 \leq Ch^{-1/4}\|\eb^\mathrm{m}_h\|_{\Lb^2}^{1/2}\left\|\jump{\curl \eb^\mathrm{m}_h}\times\nb\right\|_{\Lb^2(\Fcal_h)}^{1/2}\|\Bb\|_{\Lb^\infty}\|\eb^\mathrm{v}_I\|_{\Lb^2}.
\end{equation}
The $\Lb^2$-approximation~\eqref{eq:approxl2RT0} gives
\begin{equation}
    \begin{aligned}
        I^{\mathrm{mc}}_1 &\leq Ch^{3/4} \|\Bb\|_{\Lb^\infty} |\ub|_{\Hb^{1}} \|\eb^\mathrm{m}_h\|_{\Lb^2}^{1/2}\left\|\jump{\curl \eb^\mathrm{m}_h}\times\nb\right\|_{\Lb^2(\Fcal_h)}^{1/2}.
    \end{aligned}
\end{equation}
Since we do not assume $\Hb^1$-regularity of $\Bb$, the second term $I^{\mathrm{mc}}_2$ cannot be treated in the same way as $I^{\mathrm{mc}}_1$. Instead, we can apply the trilinear estimate~\eqref{eq:disctrilest2} since $\eb^\mathrm{m}_h \in \Zb^\mathrm{m}_h$ and $\eb^\mathrm{m}_I \in (\grad P_h)^\perp$ due to the orthogonality property~\eqref{eq:orthgradR}. Along with the interpolation approximations~\eqref{eq:approxlinfRT0},~\eqref{eq:approxl2R}, and~\eqref{eq:approxcurlR} and the inverse inequality, we obtain
\begin{equation}
    \begin{aligned}
        I^{\mathrm{mc}}_2 &= \csf(\eb^\mathrm{m}_h; \eb^\mathrm{m}_I, \ub) + \csf(\eb^\mathrm{m}_h; \eb^\mathrm{m}_I, \ub_I - \ub)\\
        &\leq C\|\ub\|_{\Wb^{1,\infty}}\|\eb^\mathrm{m}_h\|_{\Lb^2}\|\eb^\mathrm{m}_I\|_{\Hcurl} + Ch^{1/2}\|\ub\|_{\Lb^\infty}\|\eb^\mathrm{m}_I\|_{\Lb^2}\left\|\jump{\curl \eb^\mathrm{m}_h}\times\nb\right\|_{\Fcal^\circ_h} \\
        &\leq Ch^{1/2}\|\ub\|_{\Wb^{1,\infty}}(|\Bb|_{\Hb^{1/2}} + |\curl\Bb|_{\Hb^{1/2}})\|\eb^\mathrm{m}_h\|_{\Lb^2} + Ch\|\ub\|_{\Lb^\infty}|\Bb|_{\Hb^{1/2}}\left\|\jump{\curl \eb^\mathrm{m}_h}\times\nb\right\|_{\Fcal^\circ_h}.
    \end{aligned}
\end{equation}
The term $I^{\mathrm{mc}}_3$ is amenable to the trilinear estimate~\eqref{eq:disctrilest} since $\eb^\mathrm{m}_h \in \Zb^\mathrm{m}_h$. By the $\Lb^\infty$-approximation~\eqref{eq:approxlinfRT0} and the inverse inequality, we obtain
\begin{equation}\label{eq:Imcest3}
    \begin{aligned}
        |I^{\mathrm{mc}}_3| &\leq |\csf(\eb^\mathrm{m}_h; \eb^\mathrm{m}_h, \ub)| + |\csf(\eb^\mathrm{m}_h; \eb^\mathrm{m}_h, \ub_I - \ub)| \\
        &\leq C\|\ub\|_{\Wb^{1,\infty}}\|\eb^\mathrm{m}_h\|_{\Lb^2}\left(\|\eb^\mathrm{m}_h\|_{\Lb^2} + h^{1/2}\left\|\jump{\curl \eb^\mathrm{m}_h}\times\nb\right\|_{\Fcal^\circ_h}\right).
    \end{aligned}
\end{equation}
Combining the above estimates, using Young's inequality, and recalling the definition of $\ssf^\mathrm{m}_h$, we have
\begin{equation}\label{eq:Imcest}
    I^{\mathrm{mc}} \leq C_\mathrm{mc} h + C_\mathrm{mc}^\ast \|\eb^\mathrm{m}_h\|^2_{\Lb^2} + \frac{1}{4}\ssf^\mathrm{m}_h(\eb^\mathrm{m}_h, \eb^\mathrm{m}_h)
\end{equation}
where $C_\mathrm{mc}(t)$ and $C_\mathrm{mc}^\ast(t)$ depend on $\|\ub\|_{\Wb^{1,\infty}}$, $\|\Bb\|_{\Lb^\infty}$, $|\Bb|_{\Hb^{1/2}}$, and $|\curl\Bb|_{\Hb^{1/2}}$. 

\paragraph{Final estimate.}
Collecting the estimates of $I^{\mathrm{bil}}$, $I^{\mathrm{vc}}$, $I^{\mathrm{lor}}$, and $I^{\mathrm{mc}}$ in~\eqref{eq:Ibilest},~\eqref{eq:Ivcest},~\eqref{eq:Ilorest},~\eqref{eq:Imcest} and rearranging the terms, we have
\begin{multline}\label{eq:err-eq-final}
    \frac{1}{2}\frac{\mathrm{d}}{\mathrm{d}t}\left(\|\eb^\mathrm{v}_h\|^2_{\Lb^2} + \|\eb^\mathrm{m}_h\|^2_{\Lb^2}\right) + \frac{\nu_\mathrm{v}}{2}\asf_h(\eb^\mathrm{v}_h, \eb^\mathrm{v}_h) + \frac{\nu_\mathrm{m}}{2}\|\curl\eb^\mathrm{m}_h\|^2_{\Lb^2} \\+ \frac{1}{2}\ssf^\mathrm{v}_h(\eb^\mathrm{v}_h, \eb^\mathrm{v}_h) + \frac{1}{2}\ssf^\mathrm{m}_h(\eb^\mathrm{m}_h, \eb^\mathrm{m}_h) \leq C_\mathrm{err} (h + \nu_\mathrm{v} + \nu_\mathrm{m} h) + \frac{1}{2}C_\mathrm{err}^\ast\left(\|\eb^\mathrm{v}_h\|^2_{\Lb^2} + \|\eb^\mathrm{m}_h\|^2_{\Lb^2}\right)
\end{multline}
where $C_\mathrm{err}(t), C_\mathrm{err}^\ast(t) \in L^1(0, T)$ depend on $\|\partial_t \ub\|_{\Hb^{1/2+\delta}}$, $\|\ub\|_{\Wb^{1,\infty}}$, $\|\ub\|_{\Hb^{3/2+\delta}}$, $\|\partial_t \Bb\|_{\Hb^{1/2}}$, $\|\Bb\|_{\Hb^{1/2}}$, $\|\Bb\|_{\Lb^{\infty}}$, $\|\curl\Bb\|_{\Hb^{1/2+\delta}}$, and $\|\curl\Bb\|_{\Lb^\infty}$. 
The right-hand side of the above inequality is time-integrable under Assumption~\ref{ass:reg}. To simplify the notation, set 
\begin{align}
    \Ecal(t) &:= \frac{1}{2}\|\eb^\mathrm{v}_h(t)\|^2_{\Lb^2} + \frac{1}{2}\|\eb^\mathrm{m}_h(t)\|^2_{\Lb^2},\\
    \Dcal(t) &:= \frac{\nu_\mathrm{v}}{2} \asf_h(\eb^\mathrm{v}_h(t), \eb^\mathrm{v}_h(t)) + \frac{\nu_\mathrm{m}}{2} \|\curl\eb^\mathrm{m}_h(t)\|^2_{\Lb^2} + \frac{1}{2}\ssf^\mathrm{v}_h(\eb^\mathrm{v}_h(t), \eb^\mathrm{v}_h(t)) + \frac{1}{2}\ssf^\mathrm{m}_h(\eb^\mathrm{m}_h(t), \eb^\mathrm{m}_h(t)). 
\end{align}
By definition, $\Ecal(0) = 0$.
Then, the inequality~\eqref{eq:err-eq-final} simplifies to
\begin{equation}
    \frac{\mathrm{d}}{\mathrm{d}t}\Ecal(t) + \Dcal(t) \leq  (h + \nu_\mathrm{v} + \nu_\mathrm{m} h)C_\mathrm{err}(t) + \Ecal(t)C_\mathrm{err}^\ast(t).
\end{equation}
For any $t \in (0, T)$, Gronwall's inequality then gives
\begin{equation}
    \Ecal(t) + \int_0^t\exp\left(\int^t_s C_\mathrm{err}^\ast(\tau) \, \mathrm{d}\tau \right) \Dcal(s)\ds \leq (h + \nu_\mathrm{v} + \nu_\mathrm{m} h)\int_0^t\exp\left(\int^t_s C_\mathrm{err}^\ast(\tau) \, \mathrm{d}\tau \right) C_\mathrm{err}(s)\ds.
\end{equation}
Set $\overline{C}_\mathrm{err}^\ast := \int_0^T C_\mathrm{err}^\ast(s) \, \mathrm{d}s$ and $\overline{C}_\mathrm{err} := \int_0^T C_\mathrm{err}(s) \, \mathrm{d}s$. Then, we have
\begin{equation}
    \Ecal(t) + \int_0^t\Dcal(s)\ds \leq \exp\left(\overline{C}_\mathrm{err}^\ast\right)\overline{C}_\mathrm{err}(h + \nu_\mathrm{v} + \nu_\mathrm{m} h)\qquad \forall\,t \in (0, T).
\end{equation}
The triangle inequality and the $\Lb^2$-approximation~\eqref{eq:approxl2RT0},~\eqref{eq:approxl2R} then give the desired error estimate~\eqref{eq:err-est}.
\end{proof}
\begin{remark}[choice of magnetic interpolant]
    When analyzing CIP-type stabilizations (see Remark~\ref{rmk:cip}), the $\Lb^2$-orthogonal projection is normally adopted as interpolant to allow for inserting some averaging operator; see, e.g.,~\cite{burman_2006,heumann_2013,barrenechea_2024}. Here, this concerns the magnetic convection term $I^\mathrm{mc}$. We take the interpolant $\Bb_I = \Rbcal^\mathrm{ND}_h \Bb$ to exploit its good approximation in $\Hcurl$, while the role of the $\Lb^2$-projection is replaced by the discrete trilinear estimates~\eqref{eq:disctrilest} and~\eqref{eq:disctrilest2}.
\end{remark}
\begin{remark}[limitation in the diffusion-dominated regime]
    In the convection-dominated regime when $\nu_\mathrm{v} \leq h$, the estimate~\eqref{eq:err-est} implies an $O(h^{1/2})$ convergence rate. However, when the fluid diffusion is dominant, i.e., $\nu_\mathrm{v} = O(1)$, it does not imply convergence (see Figure~\ref{fig:nonaligned-convergence}).  This is expected from the poor approximation of $\Vb_h$ in the $\|\cdot\|_{\asf,h}$-norm in~\eqref{eq:asf-norm} (cf.~\cite{wang_2007}). A quick fix is to use the Brezzi--Douglas--Marini element $\Pcal_1\Lambda^2$ although the total accuracy will be limited by the magnetic discretization. Alternatively, one can resort to the vorticity-velocity formulation of the vector Laplacian~\cite{hiptmair_2026} to get a more balanced error distribution. Nevertheless, being designed for the convection-dominated regime, the proposed scheme serves well its purpose. 
\end{remark}
\begin{remark}[quasi-robustness in $\nu_\mathrm{v}$ and $\nu_\mathrm{m}$]\label{rmk:quasi-robust}
    We clarify the sense in which the error estimate~\eqref{eq:err-est} is robust with respect to $\nu_\mathrm{v}$ and $\nu_\mathrm{m}$: it is uniform over solutions satisfying Assumption~\ref{ass:reg} for which all these norms remain bounded independently of $\nu_\mathrm{v}$ and $\nu_\mathrm{m}$ (cf. the notion of \emph{quasi-robustness} in~\cite{daveiga_2025}). 
    However, it is not, by itself, $(\nu_\mathrm{v},\nu_\mathrm{m})$-uniform in general. With $\nu_\mathrm{v},\nu_\mathrm{m} \ll 1$ being perturbation parameters, the MHD system~\eqref{eq:mhd-fluid}--\eqref{eq:mhd-bc} poses a singularly-perturbed problem~\cite{temam_2018}. Boundary layers may occur and result in blow-up of solutions in high-order Sobolev norms, as the perturbation parameters approach zero~\cite{gerard_2017,liu_2019}. Nevertheless, such quasi-robust error estimates guarantee accurate approximations for smooth solutions and are commonly pursued in numerical approximations of singularly perturbed convection problems; see~\cite{roos_2008} and references therein. For parameter-uniform error estimates accommodating layer phenomena, the singular perturbation behavior of the problem must be studied and incorporated into the analysis; see, e.g.,~\cite{mardal_2002,schieweck_2008,yu_2026}. 
\end{remark}

\section{Numerical experiments}\label{sec:numexp}
In this section, we present numerical experiments to validate the theoretical results and demonstrate the performance of the proposed scheme. The implicit midpoint scheme is used for time integration and a sufficiently small time step is chosen to ensure that the temporal error is not dominant over the spatial error. The standard Picard iteration~\cite[Section~3.4]{schotzau_2004} is used to solve the nonlinear system at each time step with a relative tolerance of $10^{-10}$. The implementation is based on the finite-element library \texttt{NGSolve}~\cite{ngsolve}. The diffusion SIP penalty parameter $\gamma_\asf$ in~\eqref{eq:asfh} is set to $10$; the fluid stabilization parameter $\gamma_\mathrm{v}$ in~\eqref{eq:ssfvel} is set to $0.1$; the magnetic stabilization parameter $\gamma_\mathrm{m}$ in~\eqref{eq:ssfmag} will be specified in each test case.

We mention two slight inconsistencies with the theoretical results: the canonical interpolations~\cite[Eq.~(7.9)]{arnold_2018} are used for the discrete initial data in place of the cochain projections (see~\eqref{eq:semidiscrete}); a non-zero magnetic source term may occur in the magnetic equation~\eqref{eq:mhd-mag}. Both of these can be straightforwardly incorporated into the analysis and do not affect the convergence results.

\subsection{Case I: Smooth solutions in a convex domain}\label{sec:smoothsol-convex}

\begin{figure}[H]
\centering
\includegraphics[width=0.98\linewidth]{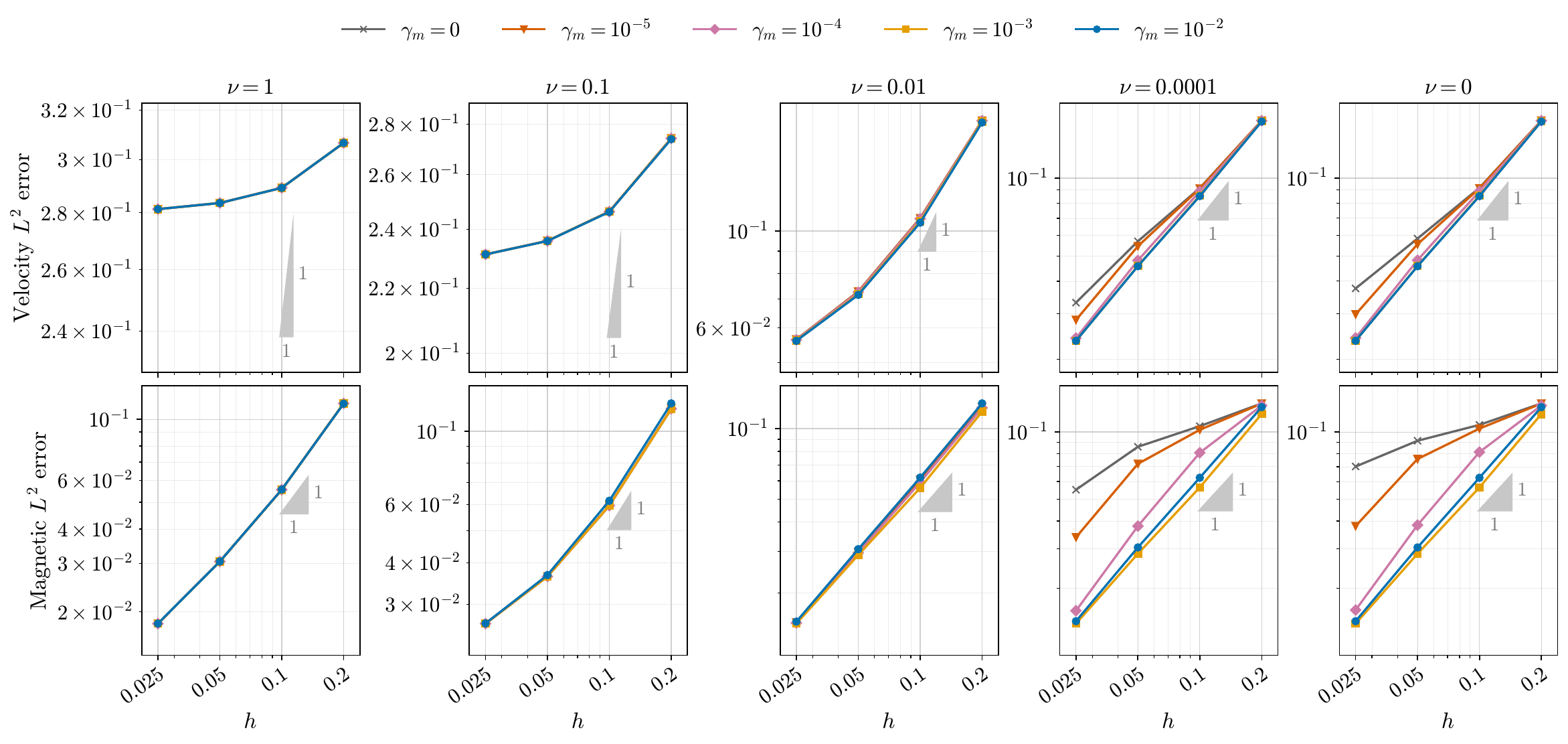}
\caption{Case I: Convergence in the velocity and magnetic $\Lb^2$-norms at final time $t = 0.1$ under mesh refinement for different diffusivity coefficients $\nu = \nu_\mathrm{v} = \nu_\mathrm{m}$ and magnetic stabilization parameter $\gamma_\mathrm{m}$.}
\label{fig:nonaligned-convergence}
\end{figure}

Let $\Omega=(0,1)^3$ be the unit cube, and the time interval $I_\mathrm{T} = (0,T)$ with $T=0.1$. Consider the manufactured solutions given by
\begin{align}
\ub(x,y,z,t)&= \frac{1 + t}{32} \curln\left(q^2(x,y,z)\cdot(1+x + 2y + 3z)\cdot[1,2,-1]^\top\right), \\
p(x,y,z,t)&= (1 + t) (x + 2 y + 3z - 3),\\
\Bb(x,y,z,t)&= \frac{1 + 2t}{48} \curln\left(q^3(x,y,z)\cdot (1 + 3x +y + 2z)\cdot[-1, 1, 2]^\top\right),
\end{align}
where $q(x,y,z) = 64x(1-x)y(1-y)z(1-z)$.
It is easy to verify that $\ub$ and $\Bb$ are divergence free and satisfy the homogeneous boundary conditions prescribed in~\eqref{eq:mhd-bc}. 
We report the final-time $\Lb^2$-errors $\|\ub(\cdot,T)-\ub_h(\cdot,T)\|_{\Lb^2}$ and $\|\Bb(\cdot,T)-\Bb_h(\cdot,T)\|_{\Lb^2}$ in Figure~\ref{fig:nonaligned-convergence} under mesh refinement for the diffusivity $\nu = \nu_\mathrm{v} = \nu_\mathrm{m} \in \{1, 10^{-1}, 10^{-2}, 10^{-4}, 0\}$ and the magnetic stabilization parameter $\gamma_\mathrm{m} \in \{0, 10^{-5}, 10^{-4}, 10^{-3}, 10^{-2}\}$. Approximately first-order convergence is observed when $\nu$ is small (i.e., $\nu = 10^{-4}, 0$). When the diffusion becomes significant (i.e., $\nu = 1, 0.1$), the velocity error saturates and the magnetic error also shows a slow-down in convergence. This is consistent with the error estimate $O(h^{1/2} + \nu_{\mathrm v}^{1/2} + \nu_{\mathrm m}^{1/2}h^{1/2})$ in~\eqref{eq:err-est} (cf.~Remark~\ref{rmk:quasi-robust}): when diffusion dominates, the poor approximation of the lowest-order Raviart--Thomas element space $\Vb_h$ in the broken $\Hb^1$-norm hinders the convergence. 

From Figure~\ref{fig:nonaligned-convergence}, it is evident that the magnetic stabilization is crucial for accuracy when $\nu$ is small. In particular, the convergence is significantly affected without the magnetic stabilization when $\nu = 10^{-4}$ and $\nu = 0$. When the diffusion is significant (i.e., $\nu = 1, 0.1$), the magnetic stabilization has a negligible effect on the convergence. 

\subsection{Case II: Singular solutions on an L-shaped domain}\label{sec:roughsol-lshape}
To test the scheme on rough solutions, we consider the benchmark problem~\cite[Test~3]{daveiga_2026b} on the L-shaped domain $\Omega=(-1,1)^3\setminus [-1,0]^2\times[-1,1]$. The manufactured stationary solutions are given by 
\begin{equation}
\ub(x,y,z) = [y^2, z^2, x^2]^{\top}, \quad p(x,y,z) = 0, \quad \Bb(x,y,z) = \grad\left(\rho^{2/3}\sin\left(\frac{2}{3}(\theta+\pi/2)\right)\right),
\end{equation}
where
$\rho=\sqrt{x^2+y^2}$ is the radial distance and $\theta\in[-\pi/2,\pi]$ is the polar angle in the $xy$-plane. The velocity $\ub$ is smooth, but the magnetic field $\Bb \in \Hb^{2/3-\delta} \;\forall\,\delta>0$ has a singularity at the reentrant corner of the L-shaped domain. 
The semi-discrete scheme~\eqref{eq:semidiscrete} is adapted to the stationary problem by replacing the time-derivative terms with the $\Lb^2$ inner products of the unknowns and the test functions. The non-zero BCs are incorporated by standard treatment. We use $\gamma_{\mathrm{m}}=0.025$ and test on $\nu_{\mathrm v}=\nu_{\mathrm m} = 10^{-6}$. Table~\ref{tab:lshape-paper-2} reports the $\Lb^2$-errors of the velocity and magnetic induction. We observe nearly optimal convergence for both fields. Note that this test case deviates from the analyzed setting in Theorem~\ref{lem:err-est} at least in two aspects: the problem is stationary and the magnetic field does not fulfill the $\Lb^\infty$ regularity in Assumption~\ref{ass:reg}. Nevertheless, the numerical results indicate that the proposed scheme works well for such rough solutions.
\begin{table}[H]
\centering\small\setlength{\tabcolsep}{3pt}
\begin{tabular}{ccccc}\hline
$h$ & $\|\ub - \ub_h\|_{\Lb^2}$ & order & $\|\Bb - \Bb_h\|_{\Lb^2}$ & order \\ \hline
0.2 & $0.26871$ & -- & $0.22616$ & --  \\
0.1 & $0.1472$ & 0.868 & $0.13365$ & 0.759 \\
0.05 & $0.077162$ & 0.932 & $0.085696$ & 0.641 \\
\hline\end{tabular}
\caption{Case II: $\Lb^2$ errors on the L-shaped domain for $\nu_\mathrm{v}=\nu_\mathrm{m}=10^{-6}$.}
\label{tab:lshape-paper-2}
\end{table}

\subsection{Case III: Periodic 2D Orszag--Tang vortex}
\begin{figure}[H]
\centering
\includegraphics[width=\linewidth]{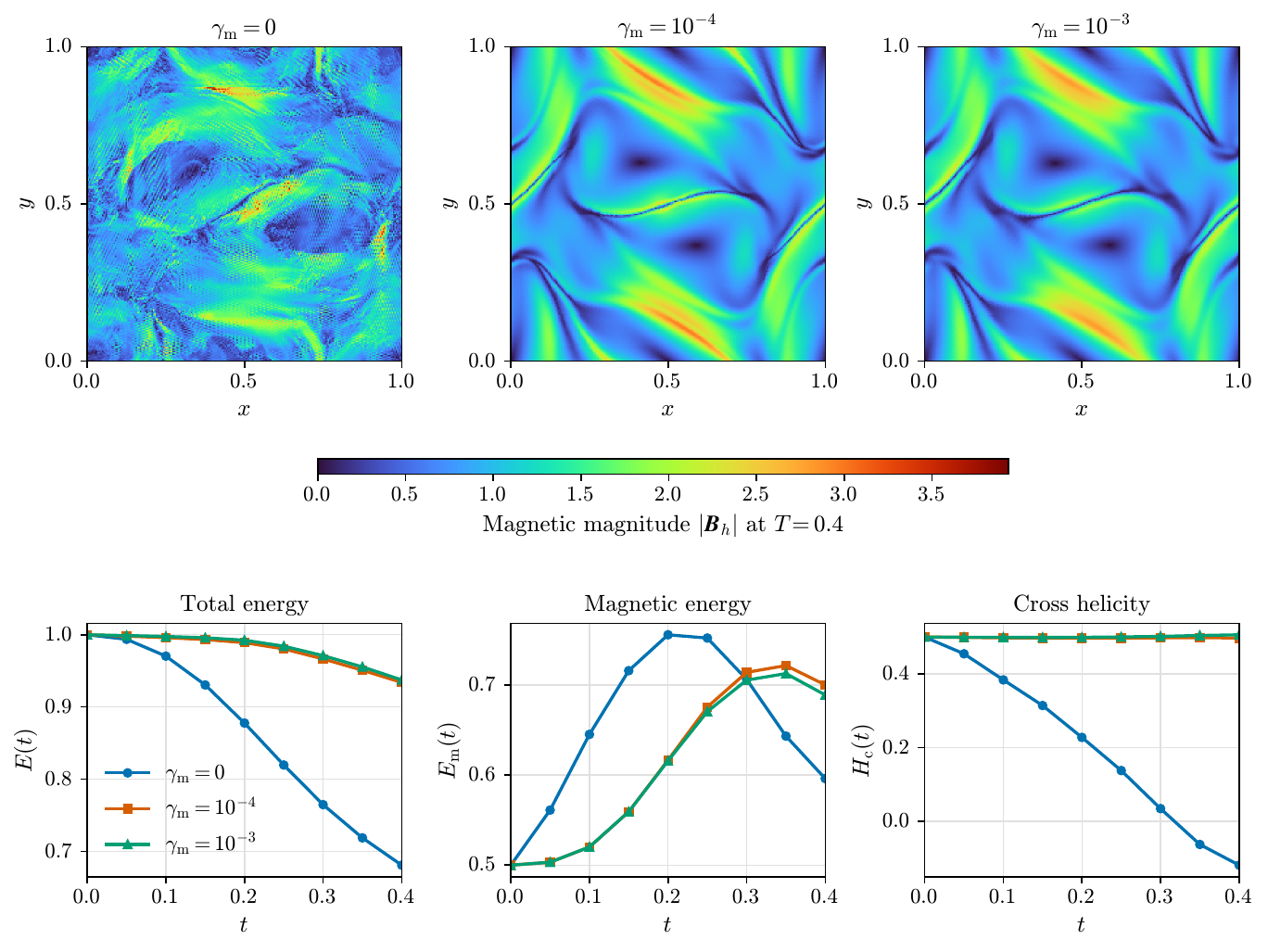}
\caption{Case III: Magnetic induction magnitude at $t=0.4$ (top) and total energy, magnetic energy, and cross helicity over time (bottom), for $\gamma_\mathrm{m} \in \{0, 10^{-4},10^{-3}\}$ for the Orszag--Tang test with mesh size $h=0.005$.}
\label{fig:orszag-tang}
\end{figure}

In this test, we consider the Orszag--Tang vortex problem~\cite{orszag_1979} on the periodic unit square $\Omega=(0,1)^2$ over the time interval $I_\mathrm{T} = (0,T)$ with $T = 0.4$ subject to the initial conditions
\begin{equation}
    \ub_0=(-\sin(2\pi y),\sin(2\pi x)), \quad
    \Bb_0=(-\sin(2\pi y),\sin(4\pi x)),
\end{equation}
subject to a zero external force.
Despite our exposition being in three dimensions with homogeneous boundary conditions, here we test the scheme~\eqref{eq:semidiscrete} in its two-dimensional form with $\curl$ reducing to $\mathrm{rot}: \Bb \mapsto \partial_x B_y - \partial_y B_x$ and with periodic boundary conditions. We fix $\nu_{\mathrm v}=\nu_{\mathrm m}=10^{-6}$ and compare different magnetic stabilization strengths $\gamma_{\mathrm m}\in \{0, 10^{-4}, 10^{-3}\}$.
The tests are conducted on a mesh with $h = 0.005$. In the top panel of Figure~\ref{fig:orszag-tang}, the magnetic induction magnitudes at $t = 0.4$ are plotted for different $\gamma_{\mathrm m}$. The lack of magnetic stabilization ($\gamma_{\mathrm m} = 0$) leads to spurious oscillations in the discrete solution. The bottom panel shows the time evolution of the total energy $E(t) := \frac{1}{2}\int_\Omega (|\ub_h|^2 + |\Bb_h|^2)\dx$, the magnetic energy $E_\mathrm{m}(t) := \frac{1}{2}\int_\Omega |\Bb_h|^2\dx$, and the cross helicity $H_\mathrm{c}(t) := \int_\Omega \ub_h\cdot\Bb_h\dx$. 
In view of the discrete energy law~\eqref{eq:energy-stability}, the total energy decay results from the velocity and magnetic diffusion as well as the introduced stabilization terms. From Figure~\ref{fig:orszag-tang}, we observe that, in the absence of the magnetic stabilization, the spurious oscillations affect the energy evolution significantly. It is known that the cross helicity is conserved in the ideal MHD case (i.e., $\nu_{\mathrm v} = \nu_{\mathrm m} = 0$)~\cite{orszag_1979}. The last figure in Figure~\ref{fig:orszag-tang} shows that the magnetic stabilization helps maintain the cross helicity, while it deviates significantly without the stabilization. Yet a discrete conservation law of the cross helicity is not addressed in this work.

\section{Concluding remarks}\label{sec:conclusion}
In this work, we present a stabilized finite element semi-discretization~\eqref{eq:semidiscrete} for the incompressible MHD equations~\eqref{eq:mhd-fluid}--\eqref{eq:mhd-bc} using the Whitney $2$-form and $1$-form elements for the velocity $\ub$ and magnetic induction $\Bb$. Apart from the exactly divergence-free velocity and weakly divergence-free magnetic induction without the aid of a magnetic Lagrange multiplier, the main contribution lies in the $(\nu_\mathrm{v}, \nu_\mathrm{m})$-independent $O(h^{1/2})$ pre-asymptotic error estimate~\eqref{eq:err-est} under the regularity assumptions on the exact solution stated in Assumption~\ref{ass:reg}. With the key choice of the CIP-type magnetic stabilization~\eqref{eq:ssfmag}, we exploit the properties of the Whitney forms (in particular, their low polynomial degrees) and the discrete de Rham complex structure. These properties allow us to show the discrete trilinear estimate~\eqref{eq:disctrilest}, which is a key step in the error analysis. The numerical experiments in Section~\ref{sec:numexp} support the convergence analysis and demonstrate the performance of the proposed scheme. Cases II and III provide additional performance tests beyond the setting of Theorem~\ref{lem:err-est}. 

This paper focuses on stabilizing the convection terms. Time integration, nonlinear solvers, and additional conservation properties remain important topics for further study.

\bibliography{ref.bib}
\end{document}